\documentclass[a4paper,12pt]{amsart}
\usepackage{amssymb,amsthm,amsmath,longtable}
\author{D.~A.~Timashev}
\address{
Department of Higher Algebra\\
Faculty of Mechanics and Mathematics\\
Moscow State University\\
119992 Moscow, Russia}
\email{timashev@mech.math.msu.su}
\title{Equivariant compactifications of reductive groups}%
\thanks{Supported by CRDF grant RM1--2088 and by RFBR grants
01--01--00756, 02--01--06402.}
\date{February 21, 2003}
\keywords{Reductive group, projective representation,
compactification, algebraic semigroup, spherical variety, orbit,
local structure, normality, smoothness}
\subjclass{14L30, 14M17, 52B20}

\makeatletter
\newcommand{\X}{\mathfrak{X}}
\newcommand{\CC}{\mathbb{C}}
\newcommand{\QQ}{\mathbb{Q}}
\newcommand{\ZZ}{\mathbb{Z}}
\newcommand{\NN}{\mathbb{N}}
\newcommand{\C}{\mathcal{C}}
\newcommand{\F}{\mathcal{F}}
\newcommand{\g}{\mathfrak{g}}
\newcommand{\lv}{\mathfrak{l}}
\newcommand{\M}{\mathcal{M}}
\renewcommand{\P}{\mathcal{P}}
\newcommand{\p}{\mathfrak{p}}
\newcommand{\V}{\mathcal{V}}
\newcommand{\D}{\mathcal{D}}
\newcommand{\PP}{\mathbb{P}}
\newcommand{\VV}{\mathbb{V}}
\newcommand{\MM}{\mathbb{M}}
\newcommand{\eps}{\varepsilon}
\newcommand{\proj}[1]{\left<{#1}\right>}
\newcommand{\supp}{\mathop{\mathrm{supp}}}
\newcommand{\intr}{\mathop{\mathrm{int}}}
\newcommand{\diag}{\mathop{\mathrm{diag}}}
\newcommand{\open}[1]{\mathaccent"7017{#1}}
\newcommand{\Hom}{\mathop{\mathrm{Hom}}}
\newcommand{\Ad}{\mathop{\mathrm{Ad}}}
\newcommand{\rk}{\mathop{\mathrm{rk}}}
\newcommand{\embeds}{\hookrightarrow}
\newcommand{\isembof}{\hookleftarrow}
\newcommand{\onto}{\twoheadrightarrow}
\newcommand{\notexists}{\m@th{\mathpalette\c@ncel\exists}}

\newcommand{\Ru}[1]{#1_{\mathrm{u}}}
\newcommand{\GL}{\mathrm{GL}}
\newcommand{\PG}{\mathrm{PG}}
\newcommand{\PGL}{\mathrm{PGL}}
\newcommand{\SL}{\mathrm{SL}}
\newcommand{\SO}{\mathrm{SO}}
\def\Sp{\mathrm{Sp}}
\newcommand{\Spin}{\mathrm{Spin}}
\newcommand{\Aa}{\mathrm{A}}
\newcommand{\Bb}{\mathrm{B}}
\newcommand{\Cc}{\mathrm{C}}
\newcommand{\Dd}{\mathrm{D}}
\newcommand{\Ee}{\mathrm{E}}
\newcommand{\Ff}{\mathrm{F}}
\newcommand{\Gg}{\mathrm{G}}
\newcommand{\LO}{\mathop{\mathrm{End}}}
\newcommand{\PO}{\mathrm{P}}
\newcommand{\Z}{\mathrm{Z}}
\newcommand{\Mat}{\mathrm{Mat}}
\newcommand{\Sym}{\mathrm{S}}
\newtheorem{theorem}{Theorem}
\newtheorem{proposition}{Proposition}
\newtheorem{lemma}{Lemma}
\newtheorem*{Corollary}{Corollary}
\theoremstyle{definition}
\newtheorem*{Definition}{Definition}
\newtheorem*{Example}{Example}
\newtheorem*{Remark}{Remark}
\makeatother

\setlongtables

\begin{document}


\begin{abstract}
We study equivariant projective compactifications of reductive
groups obtained by closing the image of a group in the space of
operators of a projective representation. We describe the
structure and the mutual position of their orbits under the action
of the doubled group by left/right multiplications, the local
structure in a neighborhood of a closed orbit, and obtain some
conditions of normality and smoothness of a compactification. Our
methods of research use the theory of equivariant embeddings of
spherical homogeneous spaces and of reductive algebraic
semigroups.

Bibliography: 36~items.
%
\end{abstract}

\maketitle

\section{Introduction}

Let $G$ be a connected reductive complex algebraic group. We may
regard $G$ as a symmetric homogeneous space $(G\times G)/\diag
G$ under the group $G\times G$ acting by left/right
multiplications. Equivariant embeddings and, in particular,
compactifications (completions) of $G$ were constructed and
studied in a number of papers. A completion of $G=\PGL_n(\CC)$
via ``complete collineations'' was constructed by
Semple~\cite{coll}. Later, Neretin generalized this construction
to other classical groups (see survey~\cite{hinges}). Neretin's
compactification may be considered as a particular case of the
``wonderful'' completion of a symmetric space, and in
particular, of a semisimple adjoint group~\cite{CP1}
(see~\S\ref{proj.comp}). A more general class of regular
completions was considered in~\cite{CP2} with applications to
the intersection theory on homogeneous spaces. The geometry and
cohomology of equivariant compactifications were studied in
\cite{cohom}, \cite{vanish}, \cite{reg.emb}, \cite{eq(wonder)},
\cite{eq(grp)}, \cite{inf(Bruhat)}, \cite{large.vars}, \cite{red.var.pro}.
The theory of equivariant embeddings of reductive groups is a
particular case of the theory of spherical varieties \cite{LV},
\cite{sph}, \cite{var.sph} (see~\S\ref{spher.emb}). On the other
hand, affine equivariant embeddings of reductive groups are
nothing else, but reductive algebraic semigroups, which were
intensively studied in a series of papers by Putcha and Renner
\cite{mon}, \cite{ss.mon}, \cite{irr.mon}, \cite{BN-mon},
\cite{ss.var}, and also by Vinberg \cite{semigr},
Rittatore \cite{rit} and Alexeev--Brion \cite{red.var.aff}.

The compactifications of Semple--Neretin and the ``wonderful''
completion of de~Concini and Procesi \cite{CP1} were based on
explicit constructions of embeddings in projective
spaces. However, the general theory of equivariant embeddings
\cite{LV} followed another way: the description of embeddings
and the study of their geometric properties goes on in terms of
certain objects of discrete convex geometry (colored cones and fans) and
of their combinatorics, see~\S\ref{spher.emb}.

Here we come back to a constructive viewpoint and consider the
following natural class of compactifications of~$G$. Let
$G:\PP(V)$ be a faithful projective representation. Then $G$
embeds in the space $\PP(\LO(V))$ of projective linear
operators, where $\LO(V)$ is the algebra of linear operators
on~$V$. The projective closure
$X=\overline{G}\subseteq\PP(\LO(V))$ is a $(G\times
G)$-equivariant projective compactification of~$G$. In fact,
$X$~is completely determined by the set of highest weights of
$G:V$.

Our goal is to extract from these data some geometric
information about $X$ including $(G\times G)$-orbit structure
(Theorem~\ref{orb}), local structure in a neighborhood of a
closed orbit (Propositions~\ref{BLV-slice},\ref{cone}),
conditions of normality (\S\ref{norm}) and smoothness
(\S\ref{smooth}). The problem of normality is important, in
particular, because we can apply a well-developed theory of
spherical varieties to normal~$X$. For instance, in this case
all orbit closures (including $X$ itself) have rational
singularities \cite{contr}, \cite{split}, there is an explicit
description of the Picard group \cite{Pic} and, in certain
cases, of the cohomology ring \cite{reg.emb}, \cite{eq(wonder)},
\cite{eq(grp)}, \cite{inf(Bruhat)}, there are vanishing theorems
for higher cohomology of nef line bundles \cite{split}, etc.

Now we explain the structure of the paper. Basic notation is
fixed in~\S\ref{not}. In \S\ref{tensor} we prove a technical
lemma on decomposing tensor powers of $G$-modules into
irreducibles. In~\S\ref{fund} we consider a particular case of
Luna's fundamental lemma, which is required in examining
smoothness of~$X$. A brief exposition of the theory of spherical
varieties, which is applied to studying projective
compactifications of~$G$, is given in~\S\ref{spher.emb}. As our
compactification $X$ is apriori a non-normal variety, it is
important to study its normalization
$\widetilde{X}\to X$. Some properties of the normalization of a
spherical variety are discussed in Proposition~\ref{biject} and
in~\S\ref{polytope} for the particular case of toric varieties.

Another important tool for the study of $X$ is the local
structure of the projective action described in \cite{BLV}. It
provides a way to construct transversal slices to closed orbits
and to reduce the study of the geometry of an action to actions
of Levi subgroups on affine subvarieties. These results are
recalled in~\S\ref{loc.struct}.

The group~$G$ as a spherical (even symmetric) homogeneous space
under $G\times G$ is considered in~\S\ref{col.data}. Here we
compute all combinatorial data required for the theory of
spherical embeddings. In \S\ref{proj.comp} we proceed to the
study of~$X$, first by describing the closed orbits and the
local structure of $X$ in their neighborhoods in terms of the
weight polytope of $G:V$. We observe that transversal slices to
closed orbits have the structure of algebraic semigroups, whose
unit groups are Levi subgroups in~$G$. These results allow to
compute combinatorial data describing~$\widetilde{X}$, which in
turn are applied to describing the orbit structure of $X$
in~\S\ref{orbits}.

Normality and smoothness of $X$ is discussed
in~\S\S\ref{norm},\ref{smooth}. Due to the local nature of
these properties, everything is reduced to the case, where $X$
is replaced by a transversal slice to a closed orbit, which is
a reductive algebraic semigroup. While an effective normality
criterion for the general case requires information on
decomposing tensor products of reductive group representations
and on branching to Levi subgroups, one can formulate some
necessary or sufficient conditions and in certain cases (e.g.\
for regular highest weights) even criteria of normality. A
criterion of smoothness is given in Theorem~\ref{reg.pt}.

In~\S\ref{examples} we illustrate the results obtained above by
the study of equivariant compactifications of simple algebraic
groups in the spaces of projective linear operators of
fundamental and adjoint representations. Some of the results
obtained here are related to similar results of Putcha--Renner
\cite{irr.mon}, \cite{ss.var} (orbital decomposition), and of
Faltings \cite{micro}, Kannan \cite{proj.norm} and de Concini
\cite{norm.sgr} (normality) for reductive algebraic semigroups.

The aim of this paper is twofold. Together with obtaining new
results, we gather in this paper and generalize some known
results on reductive group embeddings, which are scattered in
the literature. Therefore we tried to make the exposition
maximally self-contained by including proofs of some known
assertions (see e.g.~\S\ref{col.data}).

\begin{Remark}
For simplicity, we work over the field $\CC$ of complex numbers.
However our approach is purely algebraic, and all results are
valid over any algebraically closed base field of
characteristic~$0$.
\end{Remark}

The author is grateful to A.~Stoyanovsky, whose interest to the
subject stimulated the preparation of this paper, and to
G.~R\"ohrle and M.~Brion for some useful remarks, which helped to improve the
original text.

\section{Notation}
\label{not}

\noindent
$G$ is a connected reductive complex algebraic group.\\
$B\supseteq T,U$ are a fixed Borel subgroup, a maximal
torus, and a maximal unipotent subgroup in~$G$.\\
$\Delta=\Delta_G$ is the root system of $G$ relative to~$T$.\\
$\Delta^{+}\supseteq\Pi=\{\alpha_1,\dots,\alpha_l\}$ are the
subsystems of positive and simple roots relative to~$B$.\\
$\alpha^{\vee}$ is the coroot dual to~$\alpha$ (i.e., a
1-parameter subgroup in~$T$).\\
$e_{\alpha}$ is a root vector in~$\g$ corresponding to~$\alpha$.\\
$u_{\alpha}(t)=\exp(te_{\alpha})$ is a root unipotent 1-parameter subgroup.\\
$W=W_G=N_G(T)/T$ is the Weyl group.\\
$s_{\alpha}\in W$ is the reflection corresponding to
$\alpha\in\Delta$ (or its representative in~$N_G(T)$).\\
$Q=\ZZ\Delta$ is the root lattice.\\
$\X=\X(T)$ is the character lattice of~$T$.\\
$C=C_G\subseteq\X\otimes\QQ$ is the positive Weyl chamber. We
use the same notation for the positive Weyl chamber in the dual
space $\Hom(\X,\QQ)$ identified with $\X\otimes\QQ$ via a
$W$-invariant inner product.\\
$\langle\cdot,\cdot\rangle$ is the pairing between elements of
dual spaces.\\
$\X^{+}=\X^{+}_G=C\cap\X$ is the semigroup of dominant weights.\\
$\omega_1,\dots,\omega_l$ are the fundamental weights.\\
$V(\lambda)=V_G(\lambda)$ is an irreducible representation of
highest weight~$\lambda$.\\
$v_{\mu}\in V(\lambda)$ is (any) eigenvector of
$T$-weight~$\mu$. In particular, $v_{\lambda}$~is a highest
weight vector.\\
$M_{(\lambda)}$ is the isotypic component of highest
weight $\lambda$ in a $G$-module~$M$.\\
$\P=\P(V)$ is the weight polytope of a linear representation
$G:V$.\\
$\proj{v}\in\PP(V)$ is the point in the projective space
corresponding to a vector $v\in V$.\\
$\LO(V)$ is the algebra of linear operators on~$V$.\\
$\CC(X)$ is the field of rational functions on a variety~$X$.\\
$\CC[X]$ is the algebra of polynomial functions on an affine
variety~$X$.\\
$\widetilde{X}$ is the normalization of~$X$.

Lie algebras of algebraic groups are denoted by the respective
lowercase Gothic letters. The signs ``$\times,\leftthreetimes$''
denote \emph{almost} (semi)direct products of algebraic groups
(we allow intersections in finite subgroups).

For any parabolic $P\subseteq G$, $P\supseteq T$, we denote by
$P^{-}$ the opposite parabolic (and we also write $P=P^{+}$),
and by $\Ru{P}$ the unipotent radical.

\section{A result from Representation Theory}
\label{tensor}

Let $G:V$ be any rational linear representation.

\begin{lemma}[{cf.\ \cite[Lemma 4.9]{red.var.pro}}]
\label{moment}
For any $\mu\in C\cap\P$ there is $n$ such that the
decomposition of $V^{\otimes n}$ contains a simple
submodule~$V(n\mu)$.
\end{lemma}
\begin{proof}
Put $\M=\{\mu\in C\mid\exists n\in\NN:\ V^{\otimes n}\isembof
V(n\mu)\}$. Clearly, $\M\subseteq\P$. We have to prove
$\M=C\cap\P$.

Observe that $\M$ is convex. Indeed, for any $\mu,\nu\in\M$
consider any convex combination $r\mu+s\nu$, $r=p/m$, $s=q/m$,
$p,q\geq0$, $p+q=m$. For some $n$ we have $V^{\otimes
n}\isembof V(n\mu),V(n\nu)$, whence $V^{\otimes nm}\isembof
V(n\mu)^{\otimes p}\otimes V(n\nu)^{\otimes q}\isembof
V(nm(r\mu+s\nu))$. Therefore $r\mu+s\nu\in\M$.

Thus it suffices to prove $\M=C\cap\P$ for irreducible
$V=V(\lambda)$, $\lambda\in\X^{+}$. Furthermore, we may
assume $G$ to be semisimple.

First, $\lambda\in\M$. Secondly, $0\in\M$, since $V^{\otimes
n}\supset\bigwedge^nV=V(0)$ for $n=\dim V$.

Other vertices $\mu$ of $C\cap\P$ are the intersection points of
the faces of $\P$ at $\lambda$, whose direction subspaces are
spanned by simple root subsystems, with perpendicular faces
of~$C$. If $L\subset G$ is the respective Levi subgroup, then
the face is of the form
$\P_L=\P\cap(\lambda+\langle\Delta_L\rangle)$, and $\mu$ is the
center of~$\P_L$. Moreover, $V(\lambda)$~contains
$V_L(\lambda)$, whose weight polytope is~$\P_L$. Under
restriction to the commutator subgroup of $L$, $\mu$~maps to~$0$. Hence
$V(\lambda)^{\otimes n}\isembof V_L(n\mu)$. But an $L$-highest weight
vector of weight $n\mu$ in $V(\lambda)^{\otimes n}$ is
automatically a highest weight vector for~$G$ (adding
$\alpha\in\Delta^{+}\setminus\Delta_L^{+}$ moves the weight
outside the weight polytope~$n\P$), whence $\mu\in\M$.
\end{proof}

\section{Luna's fundamental lemma}
\label{fund}

The following result is an easy particular case of the
fundamental lemma \cite[II.2]{slices} in the {\'e}tale slice
theory [loc.~cit.]. For convenience of a reader, we provide a
proof.

\begin{lemma}
Suppose $X\subseteq V$ is a smooth closed $G$-stable subvariety
with a dense orbit, and $0\in X$. Then the projection
$\pi:X\to T_0X\subseteq V$ along a $G$-stable complementary
subspace is an isomorphism.
\end{lemma}
\begin{proof}
The map $\pi$ is {\'e}tale at~$0$. The set of points, where
$\pi$ is not {\'e}tale, is closed and $G$-stable. As closed
orbits are separated by polynomial invariants, $\{0\}$~is the
unique closed orbit in~$X$, hence $\pi$ is {\'e}tale. By
Zariski's Main Theorem, $\pi$~decomposes into an open immersion
$X\embeds\overline{X}$ and a finite morphism
$\overline{X}\to T_0X$, so that $\CC[\overline{X}]$ is the
integral closure of $\CC[T_0X]$ in~$\CC[X]$. But $\overline{X}$
also has a dense orbit, whence a unique closed orbit, which
means $\overline{X}\setminus X=\emptyset$. Therefore
$\pi$ is a finite {\'e}tale covering. But $\pi^{-1}(0)=\{0\}$
($X$~has a unique $G$-fixed point), hence $\deg\pi=1$,
i.e., $\pi$~is an isomorphism.
\end{proof}

\section{Spherical varieties}
\label{spher.emb}

A normal $G$-variety $X$ is called \emph{spherical} if $B$ has
an open orbit on~$X$. The more so, $X$ has an open $G$-orbit,
which can be identified with a homogeneous space $G/H$ by
choosing a base point. Thus $X$ may be considered as an
equivariant embedding of~$G/H$. All equivariant embeddings of a
spherical homogeneous space $G/H$ are described in terms of
combinatorics of certain objects from convex geometry (cones,
polytopes) related to~$G/H$. This theory is due to Luna and Vust
\cite{LV}.  For a transparent exposition, see \cite{sph},
\cite{var.sph}.

The algebra of regular functions on a spherical homogeneous
space $G/H$ is multiplicity free, i.e., nonzero isotypic
components $\CC[G/H]_{(\lambda)}\cong V(\lambda)$ are simple
$G$-modules.

Let $\X(G/H)$ be the weight lattice of all rational
$B$-eigenfunctions on~$G/H$ (which  are determined by their
weights uniquely up to proportionality). Consider the dual space
$E=\Hom(\X(G/H),\QQ)$. Any (discrete $\QQ$-valued) valuation of
$\CC(G/H)$ determines by restriction to the multiplicative
subgroup of rational $B$-eigenfunctions a homomorphism
$\X(G/H)\to\QQ$, i.e., a point in~$E$.

\begin{theorem}[{\cite[2.8, 6.3]{sph}, \cite[3.1, 4.2]{var.sph}}]
The set of $G$-invariant valuations of $\CC(G/H)$ maps to $E$
injectively, and its image is a solid polyhedral cone $\V$
containing the image of the negative Weyl chamber under the
natural projection to~$E$.
\end{theorem}

$B$-invariant divisors on~$G/H$ (i.e., irreducible components of
the complement to the open $B$-orbit, there are finitely many
of them) determine a finite set $\D$ of valuations of $\CC(G/H)$
(the set of \emph{colors}), but the map $\rho:\D\to E$ is no
more injective in general.

\begin{Definition}
A \emph{colored cone} is a pair $(\C,\F)$, where
$\F\subseteq\D$, $\rho(\F)\not\ni0$, and $\C$ is a strictly
convex polyhedral cone generated by $\rho(\F)$ and by finitely
many vectors from~$\V$, so that $\intr\C$ (the relative
interior) intersects~$\V$.

A \emph{colored face} of the colored cone is a colored cone
$(\C',\F')$, where $\C'$ is a face of $\C$ and
$\F'=\rho^{-1}(\rho(\F)\cap\C')$.

A \emph{colored fan} is a finite collection of colored cones
$(\C_i,\F_i)$ closed under passing to a colored face and such
that $\intr\C_i\cap\intr\C_j\cap\V=\emptyset$, $\forall i\ne j$
(i.e., the cones intersect in faces inside $\V$).
\end{Definition}

\begin{theorem}[{\cite[4.3]{sph}, \cite[3.4]{var.sph}}]
There is a bijection between all equivariant normal embeddings
$X\isembof G/H$ and all colored fans. Furthermore, $G$-orbits
$Y\subseteq X$ are in bijection with colored cones $(\C_Y,\F_Y)$
of a given colored fan. The set of colors $\F_Y$ corresponds to
all $B$-stable divisors on~$G/H$, whose closures contain~$Y$,
and $G$-stable divisors containing $Y$ determine the generators
of those edges of~$\C_Y$ which do not intersect $\rho(\F_Y)$.
\end{theorem}

The geometry of a spherical variety is determined by its fan.
\begin{theorem}[{\cite[4.2, 5.2]{sph}, \cite[3.4]{var.sph}}]
\label{sph.geom}
Suppose $X\isembof G/H$ is a spherical variety, and
$Y,Y'\subseteq X$ are $G$-orbits. Then
$\overline{Y}\subseteq\overline{Y'}$ iff $(\C_{Y'},\F_{Y'})$ is a
face of $(\C_Y,\F_Y)$. In particular, $G/H$~corresponds to
$(0,\emptyset)$, and $Y$ is projective iff $\C_Y$ is a solid
cone. $X$~is complete iff its fan covers all~$\V$.
\end{theorem}

Every spherical variety is
covered by \emph{simple} open subvarieties, which contain a
unique closed orbit. In particular, affine spherical varieties
are simple. The colored fan of a simple variety is determined by
the colored cone of the closed orbit and consists of all its
colored faces.

Affine spherical varieties have a more explicit description.
\begin{theorem}[{\cite[7.7]{sph}}]\label{aff}
$G/H$ admits an affine embedding iff $\rho(\D)$ is contained in
an open half-space of~$E$. There is a bijection between all
equivariant normal affine embeddings $X\isembof G/H$ and all
strictly convex polyhedral cones $\C\subset E$ generated by
$\rho(\D)$ and by finitely many vectors from~$\V$. Further,
\begin{displaymath}
\CC[X]=
\bigoplus_{\lambda\in\X(G/H)\cap\C^{\vee}}\CC[G/H]_{(\lambda)}
\end{displaymath}
and $\CC[X]^U=\CC[\X(G/H)\cap\C^{\vee}]$ is the semigroup
algebra of the semigroup of lattice points in the dual
cone~$\C^{\vee}$. The colored fan of $X$ is the set of all
colored cones of $(\C,\D)$.
\end{theorem}
\begin{Remark}
The pair $(\C,\D)$ of Theorem~\ref{aff} might not be a colored
cone in the sense of the above definition:
$\intr\C$ may have empty intersection with~$\V$.
However we may consider its colored faces in accordance with the
definition. There exists a largest face $\C'\subseteq\C$ whose
interior intersects~$\V$. The respective colored face $(\C',\F')$
corresponds to the closed orbit and determines the colored fan
of~$X$.
\end{Remark}

\begin{proposition}\label{biject}
Let $X\isembof G/H$ be a quasiprojective embedding of a
spherical homogeneous space. Then the normalization map
$\widetilde{X}\to X$ is bijective on the set of $G$-orbits.
\end{proposition}
\begin{proof}
If $X$ is affine, then $\CC[X]=
\bigoplus_{\lambda\in S}\CC[G/H]_{(\lambda)}$, where
$S\subseteq\X^{+}$ is a finitely generated semigroup, and
$\ZZ{S}=\X(G/H)$. The following objects are in bijective
correspondence: $G$-orbits on~$X$, $G$-stable closed irreducible
subvarieties of~$X$, $G$-stable prime ideals in~$\CC[X]$,
the (some, but generally not all) respective $T$-stable prime
ideals in $\CC[X]^U=\CC[S]$ (the semigroup algebra of~$S$),
(some) sets of weights of the form
$S\cap(\Sigma\setminus\Sigma')$, where $\Sigma=\QQ_{+}S$, and
$\Sigma'$ is its face.

After normalization, $S$~is replaced by $\X(G/H)\cap\Sigma$, but
the cone $\Sigma$, its faces etc.\ do not change.

The projective case is reduced to the affine case by passing to
the affine cone over $X$ and extending $G$ by homotheties. The
general case reduces to the projective case by taking the
projective closure.
\end{proof}

\begin{Example}
Toric varieties are a particular case of spherical
varieties. Here $G=B=T$ is a torus, and we may assume $H=\{e\}$,
if we want, after replacing $T$ by $T/H$. The lattice
$\X(T/H)$ coincides with the character lattice of~$T/H$. Every
$T$-invariant valuation of $\CC(T/H)$ is proportional to a
valuation given by the order in $t\to\infty$ of a function
restricted to a 1-parameter subgroup $\gamma(t)\in T/H$. Its
value at an eigenfunction of weight~$\lambda$ equals
$\langle\gamma,\lambda\rangle$, whence $\V=E$. There are no
colors, and colored cones become usual cones, and a fan is just
a finite collection of strictly convex polyhedral cones
intersecting along faces.

A simple toric variety $X$ given by a cone $\C$ is affine, and
its coordinate algebra $\CC[X]=\CC[\X(T/H)\cap\C^{\vee}]$ is the
semigroup algebra of the semigroup of lattice points in the dual
cone.

A projective toric variety can be also defined by a polytope in
$E^{*}=\X(T/H)\otimes\QQ$ which is dual to the fan, i.e., the
fan consists of dual cones to the cones at all vertices of
the polytope (and of all their faces).
\end{Example}

\section{The polytope of the closure of a torus orbit}
\label{polytope}

Consider a linear representation $T:V$ and a vector $v\in V$.
The \emph{support} $\supp{v}$ is the convex hull of all weights
in the weight decomposition of~$v$.

\begin{proposition}
The polytope of the projective toric variety
$\widetilde{\overline{T\proj{v}}}$ equals $-\supp{v}$.
\end{proposition}
\begin{proof}
Consider the weight decomposition
$v=v_{\lambda_0}+\dots+v_{\lambda_n}$, complete the set
$v_{\lambda_0},\dots,v_{\lambda_n}$ to a weight basis of~$V$,
and consider the dual coordinates $x_{\lambda}$ (of
weights~$-\lambda$) as homogeneous coordinates in~$\PP(V)$. The
variety $X=\overline{T\proj{v}}$ is covered by affine charts
$X_i=\{x_{\lambda_i}\ne0\}$. To be definite, consider~$X_0$.
Then $\CC[X_0]=
\CC[x_{\lambda_1}/x_{\lambda_0},\dots,x_{\lambda_n}/x_{\lambda_0}]$
is the semigroup algebra of the semigroup
$S_0=\ZZ_{+}(\lambda_0-\lambda_1)+\dots+\ZZ_{+}(\lambda_0-\lambda_n)$.
This semigroup generates the cone $\Sigma_0$ at the vertex
$-\lambda_0$ of~$-\supp v$. Hence
$\CC[\widetilde{X_0}]=\CC[\QQ{S_0}\cap\Sigma_0]$. Thus the fan
of $\widetilde{X}$ is dual to~$-\supp{v}$.
\end{proof}

\begin{Corollary}[of the proof]
$\overline{T\proj{v}}$ is normal iff
$\QQ{S_i}\cap\Sigma_i=S_i$, $i=0,\dots,n$
\end{Corollary}

\section{The local structure of a projective action}
\label{loc.struct}

Consider a linear representation $G:V$ and fix a lowest weight vector
$v_{\lambda}\in V$. Let
$P^{-}=P(\lambda)=G_{\proj{v_{\lambda}}}$, and $P\supseteq B$ be
the opposite parabolic. Consider the Levi decompositions
$P=\Ru{P}\leftthreetimes L$,
$P^{-}=\Ru{P^{-}}\leftthreetimes L$. We have a decomposition
$V=\Ru{\p}v_{\lambda}\oplus\CC{v_{\lambda}}\oplus M$ in a direct
sum of $L$-modules. Consider the dual highest weight vector
$v_{-\lambda}\in V^{*}$ and the open subset
$\open{V}=\Ru{\p}v_{\lambda}\oplus\CC^{\times}{v_{\lambda}}\oplus M$,
which is the complement to the hyperplane
$\{\langle v_{-\lambda},\cdot\rangle=0\}$.

\begin{theorem}[\cite{BLV}]
The $P$-action on $\CC{v_{\lambda}}\oplus M$ yields an
isomorphism
\begin{displaymath}
\PP(\open{V})\cong
P\times_{L}M(-\lambda)=\Ru{P}\times M(-\lambda)
\end{displaymath}
where $M(-\lambda)=M$ is equipped with the $L$-action twisted by
the character $-\lambda$ and embedded in $\PP(V)$ as
$\{\proj{v_{\lambda}}\}\times M=\PP(v_{\lambda}+M)$.
\end{theorem}
It follows that any projective $G$-variety $X\subseteq\PP(V)$
containing the closed orbit $Y=G\proj{v_{\lambda}}$ contains the
affine open $P$-stable subset $\open{X}=X\cap\PP(\open{V})\cong
P\times_{L}Z=\Ru{P}\times Z$ intersecting $Y$, where $Z\subseteq
M(-\lambda)$ is a closed affine $L$-subvariety.

If $X$ is a spherical $G$-variety, then $Z$ is an affine
spherical $L$-variety. It is easily seen from the structure
of~$\open{X}$ that $\X(X)=\X(\open{X})=\X(Z)$, and
$\C_Y=\C_{Y\cap Z}$. The only delicate point is that the cone of
invariant valuations can increase and the set of colors can
decrease, because some colors may become $L$-stable divisors.
Combined with Theorem~\ref{aff}, this yields an effective
description of the local structure of a spherical variety by its
fan.

\section{A reductive group as a spherical homogeneous space}
\label{col.data}

The group $G$ is a homogeneous space under $G\times G$ acting by
left and right multiplications with stabilizer $\diag{G}$
of~$e$. Fix a Borel subgroup $B^{-}\times B\subseteq
G\times G$. The Bruhat decomposition implies that $G=(G\times
G)/\diag{G}$ is a spherical homogeneous space.
The combinatorial data related to this space in the sense of
\S\ref{spher.emb} were computed by Vust \cite{symm} (in the more
general context of symmetric spaces) and by Rittatore \cite{rit}.
We reproduce their results below. The following result is well
known \cite[II.3.1, Satz~3]{inv}:

\begin{proposition}
$\CC[G]=\bigoplus_{\lambda\in\X^{+}}\CC[G]_{(\lambda)}$, where
$\CC[G]_{(\lambda)}\cong V(\lambda)^{*}\otimes V(\lambda)$ is the
linear span of the matrix elements of the representation
$G:V(\lambda)$.
\end{proposition}

\begin{Corollary}
$\X(G)\cong\X$
\end{Corollary}

The eigenfunction $f_{\lambda}(g)=\langle
v_{-\lambda},gv_{\lambda}\rangle$ has highest weight
$(-\lambda,\lambda)$, i.e., $\lambda$ under our identifications.

\begin{proposition}\label{tails}
$\CC[G]_{(\lambda)}\cdot\CC[G]_{(\mu)}=\CC[G]_{(\lambda+\mu)}
\oplus\bigoplus_i\CC[G]_{(\lambda+\mu-\beta_i)}$, where
$\lambda+\mu-\beta_i$ are highest weights of all ``lower''
irreducible components in the decomposition $V(\lambda)\otimes
V(\mu)=V(\lambda+\mu)\oplus\dots$, so that
$\beta_i\in\ZZ_{+}\Pi$.
\end{proposition}
\begin{proof}
$\CC[G]_{(\lambda)}\cdot\CC[G]_{(\mu)}$ is generated by
products of matrix elements of $G:V(\lambda)$ and $G:V(\mu)$,
i.e., by matrix elements of $G:V(\lambda)\otimes V(\mu)$.
\end{proof}

\begin{Corollary}
$\V=-C$ is the negative Weyl chamber.
\end{Corollary}
\begin{proof}
It is easy to see that a $G$-invariant valuation of $\CC(G)$ is
constant at any isotypic component $\CC[G]_{(\lambda)}$, and its
restriction is a linear function of $\lambda$, i.e.,
$\nu\in E=\Hom(\X,\QQ)$. The value of the valuation at $\forall
f\in\CC[G]$ equals
$\min_{f_{(\lambda)}\ne0}\langle\nu,\lambda\rangle$, where
$f_{(\lambda)}$~is the projection of~$f$
on~$\CC[G]_{(\lambda)}$. Multiplying functions from
$\CC[G]_{(\lambda)}$ and $\CC[G]_{(\mu)}$, we deduce from
Proposition~\ref{tails} that $\langle\nu,\beta_i\rangle\leq0$ is
necessary for $\nu$ to define a valuation. Since
$\beta_i$ are positive combinations of simple roots, and
multiples of all simple roots are among them for appropriate
$\lambda,\mu\in\X^{+}$ (see e.g.~\S\ref{tensor}), these
inequalities define~$-C$.

Conversely, each $\nu\in-C$ defines a $\QQ$-valued function on
$\CC[G]$ (denoted by the same letter by abuse of notation)
satisfying the additive property of a valuation by the above
formula. To verify the multiplicative property, take
$p,q\in\CC[G]$ and choose $\gamma\in\intr(-C)$ such that
$\min\langle\gamma,\lambda\rangle$ over all $\lambda$ with
$p_{(\lambda)}\ne0$, $\langle\nu,\lambda\rangle=\min$, and
$\min\langle\gamma,\mu\rangle$ over all $\mu$ with
$q_{(\mu)}\ne0$, $\langle\nu,\mu\rangle=\min$, are reached at the
unique points $\lambda_0$ and $\mu_0$, respectively. Then
$pq=f_{(\lambda_0+\mu_0)}+\sum f_{(\chi)}$, where
$\langle\nu,\chi\rangle\ge\langle\nu,\lambda_0+\mu_0\rangle$,
$\langle\gamma,\chi\rangle\ge\langle\gamma,\lambda_0+\mu_0\rangle$,
and at least one of these inequalities is strict. Hence
$\nu(pq)=\langle\nu,\lambda_0+\mu_0\rangle=\nu(p)+\nu(q)$, and we
are done.
\end{proof}

It follows from the Bruhat decomposition, that $(B^{-}\times
B)$-stable divisors on~$G$ are of the form
$D_i=\overline{B^{-}s_{\alpha_i}B}$. If $G$ is of simply
connected type (i.e., $G$~is a direct product of a torus and a
simply connected semisimple group), then $\omega_i\in\X$, and
$D_i$ is defined by the equation $f_{\omega_i}=0$. More
precisely, consider a curve
$g_j(t)=u_{-\alpha_j}(t)\alpha_j^{\vee}(t^{-1})u_{\alpha_j}(-t)$
in the big Bruhat cell. Then
$\lim_{t\to\infty}g_j(t)=s_{\alpha_j}\in D_j$ (everything takes
place in an $\SL_2$-subgroup, where it is computed explicitly),
and
\begin{displaymath}
f_{\omega_i}(g_j(t))=
\langle v_{-\omega_i},\alpha_j^{\vee}(t^{-1})v_{\omega_i}\rangle=
t^{-\langle\omega_i,\alpha_j^{\vee}\rangle}
\langle v_{-\omega_i},v_{\omega_i}\rangle
\end{displaymath}
Hence $f_{\omega_i}$ has order $1$ along $D_j$ for $i=j$ and
$0$ for $i\ne j$. Therefore the colors look like
$\rho(D_j)=\alpha^{\vee}_j$. This conclusion remains valid for
arbitrary~$G$, because $G$ is covered by a group of simply
connected type with the colors being the preimages of the colors
of~$G$.

These results, due to Vust~\cite{symm}, allow us to apply the
theory of \S\S\ref{spher.emb},\ref{loc.struct} to describing
$(G\times G)$-equivariant embeddings of~$G$. Affine embeddings
have another remarkable property.

\begin{theorem}[\cite{semigr}]
An affine $(G\times G)$-equivariant embedding $X\isembof G$ is
an algebraic semigroup with unit, and $G$ is the group of
invertibles in~$X$.
\end{theorem}
\begin{proof}
The actions of the left and right copy of $G\times G$ on $X$
define coactions $\CC[X]\to\CC[G]\otimes\CC[X]$ and
$\CC[X]\to\CC[X]\otimes\CC[G]$, which are the restrictions to
$\CC[X]\subseteq\CC[G]$ of the comultiplication
$\CC[G]\to\CC[G]\otimes\CC[G]$. Hence the image of $\CC[X]$ lies
in
$\CC[G]\otimes\CC[X]\cap\CC[X]\otimes\CC[G]=\CC[X]\otimes\CC[X]$,
and we have a comultiplication in~$\CC[X]$. Now $G$ is open in
$X$ and consists of invertibles. For any invertible $x\in X$, we
have $xG\cap G\ne\emptyset$ $\implies x\in G$.
\end{proof}

\section{Projective compactification of a reductive group}
\label{proj.comp}

Let $V=\bigoplus_{i=0}^{m}V(\lambda_i)$
($\lambda_i\in\X^{+}$) be a faithful representation of~$G$. Then
$G\embeds\VV=\bigoplus\LO(V(\lambda_i))\subseteq\LO(V)\cong
V\otimes V^{*}$. Passing to the projectivization, we have
$G\embeds\PP=\PP(\VV)\iff\X=\sum\ZZ(\lambda_i-\lambda_j)+Q$.
This condition can always be achieved by adding to $V$ a trivial
representation, and we will assume it or, more generally, we
will consider faithful projective representations $G:\PP(V)$
(which can be lifted to a linear representation of rather a
finite cover of~$G$, than $G$~itself). Moreover, we may assume
all $\lambda_i$ to be distinct.

Our aim is to study $X=\overline{G}\subseteq\PP$.

The image of the identity in $\VV$ is $e=\sum v_{\mu}\otimes
v_{-\mu}$, where $v_{\pm\mu}$ run over dual weight bases in
$V,V^{*}$. Under the above identification $\X((T\times
T)/\diag{T})\cong\X(T)$, we have $\supp{e}=-\P$.

All closed $(G\times G)$-orbits in $\PP$ are of the form
$Y_i=(G\times G)y_i$, $y_i=\proj{v_{\lambda_i}\otimes
v_{-\lambda_i}}$.

The local structure of $\PP$ in a neighborhood of $Y_i$ looks as
follows (\S\ref{loc.struct}). Let $P=P(\lambda_i)$ be the
projective stabilizer of $v_{\lambda_i}$, and
$P=\Ru{P}\leftthreetimes L$, $P^{-}=\Ru{P^{-}}\leftthreetimes L$
be the Levi decompositions. Then $P\times P^{-}=(G\times
G)_{y_i}$. Consider the open subset
\begin{displaymath}
\open{\VV}=\CC^{\times}(v_{\lambda_i}\otimes v_{-\lambda_i})
\oplus(\Ru{\p^{-}}v_{\lambda_i}\otimes v_{-\lambda_i})\oplus
(v_{\lambda_i}\otimes\Ru{\p}v_{-\lambda_i})\oplus\MM
\end{displaymath}
in~$\VV$, where $\MM\subset\VV$ is an $(L\times L)$-stable
subspace and $\CC(v_{\lambda_i}\otimes v_{-\lambda_i})\oplus\MM
\ni{e}$. The affine $(P^{-}\times P)$-stable chart
$\open{\PP}=\PP(\open{\VV})$ intersects $Y_i$, so that
$\open{\PP}\cong(P^{-}\times P)\times_{(L\times
L)}\MM(-\lambda_i,\lambda_i)=
\Ru{P^{-}}\times\Ru{P}\times\MM(-\lambda_i,\lambda_i)$

\begin{proposition}
$X\supset Y_i\iff\lambda_i$ is a vertex of~$\P$
\end{proposition}
\begin{proof}
If $\lambda_i$ is a vertex of~$\P$, then there exists a
1-parameter subgroup $\gamma(t)\in T$ such that
$\langle\gamma,\lambda_i\rangle>\langle\gamma,\mu\rangle$,
Then
\begin{displaymath}
\gamma(t)=\proj{\sum\gamma(t)v_{\mu}\otimes v_{-\mu}}=\proj{\sum
t^{\langle\gamma,\mu\rangle}v_{\mu}\otimes v_{-\mu}}\to y_i
\qquad\mbox{as }t\to\infty
\end{displaymath}
Conversely, assume $Y_i\subset X$. Consider the local structure
of $X$ in a neighborhood of~$Y_i$. $X$~contains the
$(P^{-}\times P)$-stable affine open subset
$\open{X}=\Ru{P^{-}}\times\Ru{P}\times Z$, where
$Z=X\cap{\MM(-\lambda_i,\lambda_i)}=\overline{(L\times L)e}
=\overline{L}$ is an affine $(L\times L)$-equivariant embedding
of $L$ with the unique fixed point~$y_i$. Then by
Theorem~\ref{aff}, $\CC[\widetilde{Z}]=
\bigoplus_{\lambda\in\X\cap\C^{\vee}}\CC[L]_{(\lambda)}$,
where $\C$ is a strictly convex solid cone generated by
$\Delta^{+}_L$ and by finitely many vectors from $\V_L=-C_L$,
so that $\intr\C\cap\V_L\ne\emptyset$. Hence
$\CC[Z]=\bigoplus_{\lambda\in S}\CC[L]_{(\lambda)}$, where
$S\subseteq\X^{+}_L\cap\C^{\vee}$ is a finitely generated
semigroup such that $\ZZ{S}\cap\C^{\vee}=\X\cap\C^{\vee}$.

There exists a 1-parameter subgroup
$\gamma\in\intr\C\cap-C_L$, $\gamma\perp\Delta_L$. It defines a
non-negative grading of~$\CC[Z]$:
$\deg\CC[Z]_{(\lambda)}=\langle\gamma,\lambda\rangle$, and
functions of zero degree are constant. Consequently the action
of $\gamma(t)$ contracts the whole $Z$ (and in particular,~$e$)
to $y_i$ as $t\to\infty$, whence $-\lambda_i$ is a vertex of
$\supp{e}=-\P$.
\end{proof}

Let us give a more detailed study of the local structure of $X$
in a neighborhood of a closed orbit, say~$Y_0$. We keep the
previous notation. Put
$\open{X}=X\cap\open{\PP}=\{f_{\lambda_0}\ne0\}$.
(Here $f_{\lambda_0}=v_{-\lambda_0}\otimes v_{\lambda_0}$
is regarded as a linear function on~$\VV$.) Then
$\open{X}\cong(P^{-}\times P)\times_{(L\times
L)}Z=\Ru{P^{-}}\times\Ru{P}\times Z$, where
$Z=X\cap(\{y_0\}\times\MM)=\overline{L}$ is an equivariant
affine embedding.

More precisely, if $V=\CC v_{\lambda_0}\oplus V_0$ is an
$L$-stable decomposition, then $Z\subseteq\LO(V_0(-\lambda_0))
\subseteq\MM(-\lambda_0,\lambda_0)$, whence
\begin{proposition}\label{BLV-slice}
Let $\mu_1,\dots,\mu_s\in\X^{+}_L$ be all highest weights
of $L:V$, except~$\lambda_0$. The representations
$L:V_L(\mu_j-\lambda_0)$ extend to $Z$ and induce a closed
embedding
$Z\embeds\VV_0=\bigoplus\LO(V_L(\mu_j-\lambda_0))$.
\end{proposition}

\begin{proposition}\label{cone}
$\CC[Z]=\bigoplus_{\lambda\in S}\CC[L]_{(\lambda)}$, where
$S\subseteq\X^{+}_L$ is a semigroup generating the cone
$\Sigma_0$ of~$C\cap\P$ at~$\lambda_0$, and $\ZZ{S}=\X$.
\end{proposition}
\begin{proof}
The homogeneous coordinate algebra of~$\PP$ is
$\CC[\VV]=\Sym^{\bullet}\VV^{*}$. The set of weights of
homogeneous polynomials of degree $n$ is contained in
$-n\P\oplus n\P$. As
$\CC[\open{\PP}]=\bigoplus\Sym^n\VV^{*}/f_{\lambda_0}^n$, the
weights of homogeneous polynomials of degree $n$ on
$\open{\PP}$ lie in $n(-\P+\lambda_0)\oplus n(\P-\lambda_0)$.
Hence the highest weights of the $(L\times L)$-module
$\CC[Z]\subseteq\CC[L]$ lie in $C_L\cap\bigcup
n(\P-\lambda_0)=\Sigma_0$ (under the antidiagonal embedding
$\X\subset\X(T\times T)$).

On the other hand, for each $\mu\in C_L\cap\P$,
Lemma~\ref{moment} yields $V^{\otimes n}\isembof V_L(n\mu)$ for
some~$n$. Using Proposition~\ref{tails} for~$L$, it is easy to
derive $f_{n\mu}/f_{\lambda_0}^n\in\CC[Z]$. Therefore $S$
generates $\Sigma_0$.
\end{proof}

\begin{Corollary}
$\CC[\widetilde{Z}]=
\bigoplus_{\lambda\in\X\cap\Sigma_0}\CC[L]_{(\lambda)}$
\end{Corollary}

\begin{Remark}
If $\lambda_0$ is a regular weight (i.e., $\lambda_0\in\intr
C$), then $L=T$, $\Sigma_0$ is the cone of $\P$ at~$\lambda_0$,
and $\CC[Z]=\CC[S]$ is the semigroup algebra of the semigroup
$S$ generated by $\mu-\lambda_0$, where $\mu$ runs over all
weights of~$V$ (cf.~\S\ref{polytope}). In the general case, the
description of~$S$ requires decomposing the $G$-module $V$
into simple $L$-modules and also tensor products of simple
$L$-modules (see~\S\ref{norm}). Normality of $Z$ means
$S=\X\cap\Sigma_0$.
\end{Remark}

Consider the spherical variety $\widetilde{X}$. By
Proposition~\ref{biject}, the normalization map
$\widetilde{X}\to X$ is bijective on the set of orbits. In
particular, $\widetilde{X}\isembof G/H,Y_i$.

\begin{theorem}
\label{fan}
The colored fan of $\widetilde{X}$ is generated by colored cones
of closed orbits $Y_i\subset\widetilde{X}$ over all vertices
$\lambda_i\in\P$: $\C_{Y_i}=\Sigma_i^{\vee}$, where
$\Sigma_i$ is the cone of~$C\cap\P$ at~$\lambda_i$, and
$\F_{Y_i}=\Pi_{L_i}^{\vee}$ is the subsystem of simple coroots
of the respective Levi subgroup $L_i=Z_G(\lambda_i)$. (Recall
that colors are identified with simple coroots.)
\end{theorem}
\begin{proof}
The assertion on $\C_{Y_i}$ stems from Proposition~\ref{cone}
and from arguments in~\S\ref{loc.struct}. Let us deal with
colors.

Note that $D_j\in\F_{Y_i}\iff y_i\in\overline{D_j}$. As
$f_{\lambda_i}(y_i)\ne0$ and $f_{\lambda_i}$ has order
$\langle\lambda_i,\alpha^{\vee}_j\rangle$ along $D_j$, it
follows that $D_j\notin\F_{Y_i}$ whenever
$\lambda_i\not\perp\alpha_j\iff\alpha_j\notin\Pi_{L_i}$.
Conversely, assume $\alpha_j\in\Pi_{L_i}$, and let $\gamma(t)\in
T$ be a 1-parameter subgroup such that
$\langle\gamma,\lambda_i\rangle>\langle\gamma,\mu\rangle$,
$\forall\mu\in\P$, $\mu\ne\lambda_i$. Then as $t\to\infty$,
\begin{displaymath}
\gamma(t)s_{\alpha_j}=
\proj{\sum\gamma(t)v_{s_{\alpha_j}\mu}\otimes v_{-\mu}}=
\proj{\sum t^{\langle\gamma,s_{\alpha_j}\mu\rangle}
v_{s_{\alpha_j}\mu}\otimes v_{-\mu}}\to y_i,
\end{displaymath}
i.e., $\overline{D_j}\ni y_i$.
\end{proof}

\begin{Remark}
In case, where all vertices of $\P$ are regular weights, the
fan of $X$ has no colors and is a subdivision of~$-C$.
Transversal slices $Z_i$ to closed orbits $Y_i$ at $y_i$ are
affine toric varieties, and their cones are dual to the cones of
$\P$ at $\lambda_i$. These slices are contained in a projective
toric variety $\overline{T}\subseteq X$ (or $\widetilde{X}$).
Clearly, $\overline{T}$ is $W$-stable ($W$~or, more precisely,
$N_G(T)$ acts by conjugation). It follows from the local
structure that $\overline{T}$ intersects all $Y_i$
transversally at the points $w(y_i)$, $w\in W$, which are all
$T$-fixed points of~$\overline{T}$. The fan of~$\overline{T}$ is
the $W$-span of the fan of~$X$, and its polytope equals $\P$ by
Proposition~\ref{polytope}.

Left cosets $\overline{vT}$, $v\in W$, intersect all $Y_i$
transversally at the points $(vw,w)y_i$, $w\in W$, which are all
$(T\times T)$-fixed points of~$X$. The closure of the normalizer
of the torus $\overline{N_G(T)}=\bigsqcup_{v\in W}\overline{vT}$
is a disjoint union of cosets. (Otherwise distinct
$\overline{vT}$ would share fixed points.)

In case, where $V=V(\lambda_0)$ is irreducible, $X$~is normal
(even smooth, cf.~\S\S\ref{norm},\ref{smooth} and \cite{CP1}) and does not
depend on the choice of (regular) $\lambda_0$. This
compactification is called \emph{wonderful}~\cite{CP1}. The
general case of regular weights was considered in \cite{CP2}
and in a number of other papers.
\end{Remark}

\section{Orbits}
\label{orbits}

According to the general theory of spherical varieties and
Theorem~\ref{fan}, $(G\times G)$-orbits in $X$ whose closures
contain a given closed orbit $Y_i$ are in bijection with colored
faces of the colored cone $(\Sigma_i^{\vee},\Pi_{L_i}^{\vee})$
or, equivalently, with faces $\C\subseteq\Sigma_i^{\vee}$ such
that $(\intr\C)\cap(-C)\ne\emptyset$. We reformulate this
description in ``dual'' terms of~$\P$.

For any face $\Gamma\subseteq\P$, denote by $|\Gamma|$ its
direction subspace, by $\langle\Gamma\rangle$ its linear span,
and put
$\|\Gamma\|=|\Gamma|\oplus(\langle\Gamma\rangle^{\perp}\cap\langle
Q\rangle)$. (If $G$ is semisimple and $\Gamma$ is a proper face,
then $\|\Gamma\|$ is its supporting hyperplane shifted to~$0$.) We
say that $\gamma\in\Hom(\X,\QQ)$ is a \emph{supporting function}
for $\Gamma$ in $\P$ if $\langle\gamma,\Gamma\rangle=\text{const}<
\langle\gamma,\P\setminus\Gamma\rangle$, i.e., $\gamma$~is the
l.h.s.\ of the equation of a supporting hyperplane for $\Gamma$.

\begin{proposition}\label{orb<->faces}
There is a bijection between all $(G\times G)$-orbits
$Y\subset X$ and all faces $\Gamma\subseteq\P$ such that $(\intr\Gamma)\cap
C\ne\emptyset$. Here $\C_Y$ is dual to the cone of
$C\cap\P$ at the face $C\cap\Gamma$, and
$\F_Y=\{D_j\mid\alpha_j\perp\langle\Gamma\rangle\}$. ``Colored''
orbits (i.e., those having $\F_Y\ne\emptyset$) correspond to
faces lying at the boundary of~$C$. The adherence of orbits
(i.e., inclusion of orbit closures) corresponds to the inclusion
of the respective faces.
\end{proposition}
\begin{proof}
Assume $\overline{Y}\supseteq Y_i$, i.e., $\C_Y$~is a face
of~$\Sigma_i^{\vee}$. We have $\gamma\in\intr\C_Y\iff$
$\gamma$~is a supporting function for the dual face
$\Sigma_Y=\Sigma_i\cap\C_Y^{\perp}$ of~$\Sigma_i$, i.e.,
$\langle\gamma,\Sigma_i\rangle\geq0$ and $\Sigma_Y$ is
distinguished by the equation $\langle\gamma,\cdot\rangle=0$.
Thus colored faces $\Sigma_i^{\vee}$ correspond to faces
of~$\Sigma_i$, or of $C\cap\P$ at~$\lambda_i$, cut out by
supporting functions $\gamma\in-C$. But such faces of $C\cap\P$ are
obtained by intersecting $C$ with faces $\Gamma\subseteq\P$ of
the same dimension cut out by $\gamma$ in~$\P$.

The description of colors stems from
$\Pi_{L_i}^{\vee}=\Pi^{\vee}\cap\lambda_i^{\perp}$ and
$\F_Y=\Pi_{L_i}^{\vee}\cap|\Gamma|^{\perp}$. The assertion on
adherence follows from Theorem~\ref{sph.geom}, because duality
reverts inclusion of faces.
\end{proof}

Let us make the above description of orbits more explicit. We
introduce the following notation.

For any subspace $N\subseteq\X\otimes\QQ$ orthogonal to some
dominant weight, denote by $P_N$ the parabolic in~$G$ generated
by $B$ and by the roots $\alpha\in N$. Let
$P_N^{\pm}=\Ru{(P_N^{\pm})}\leftthreetimes L_N$ be the Levi
decompositions, $L_N\supseteq T$, $L_N'$~be the commutator
subgroup of~$L_N$. For any sublattice $\Lambda\subseteq\X$, denote by
$T^{\Lambda}\subseteq T$ the diagonalizable group which is the common
kernel of all characters $\lambda\in\Lambda$.

For any face $\Gamma\subseteq\P$, denote by $V_{\Gamma}$ the sum
of weight subspaces of~$V$ with weights in~$\Gamma$, and by
$V_{\Gamma}'$ its $T$-stable complement. Let
$e_{\Gamma}=\sum_{\mu\in\Gamma}v_{\mu}\otimes v_{-\mu}$ be the
projector on $V_{\Gamma}$ w.r.t.\ the decomposition
$V=V_{\Gamma}\oplus V_{\Gamma}'$.

Observe that $V_{\Gamma},V_{\Gamma}'$ are $L_{|\Gamma|}$- and
even $L_{\|\Gamma\|}$-stable,
$L_{\|\Gamma\|}'=L_{|\Gamma|}'\times
L_{\langle\Gamma\rangle^{\perp}}'$, and the action
$L_{\langle\Gamma\rangle^{\perp}}':V_{\Gamma}$ is trivial.
Indeed, adding roots
$\alpha\in\|\Gamma\|\setminus|\Gamma|$ moves $\Gamma$
outside~$\P$, i.e., the respective root vectors act on
$V_{\Gamma}$ trivially. This means that
$\Delta\cap\|\Gamma\|=(\Delta\cap|\Gamma|)\sqcup
(\Delta\cap\langle\Gamma\rangle^{\perp})$ is a disjoint
orthogonal union.

In the subspaces $|\Gamma|,\langle\Gamma\rangle$, consider
generating sublattices
\begin{align*}
|\Gamma|_{\ZZ}&=
\sum_{\lambda_i,\lambda_j\in\Gamma}\ZZ(\lambda_i-\lambda_j)+
Q\cap|\Gamma|\\
\langle\Gamma\rangle_{\ZZ}&=
\sum_{\lambda_i\in\Gamma}\ZZ\lambda_i+Q\cap|\Gamma|
\end{align*}
(To be rigorous, these lattices depend not only on $\P$
and~$\Gamma$, but on the initial set of highest weights
$\lambda_0,\dots,\lambda_m$.) Note that
$\langle\Gamma\rangle_{\ZZ}$ is the weight lattice of
$T:V_{\Gamma}$.

\begin{theorem}\label{orb}
The orbit $Y\subset X$ corresponding to the face
$\Gamma\subseteq\P$ is represented by $y=\proj{e_{\Gamma}}$.
Stabilizers look like
\begin{align*}
(G\times G)_y&=
\bigl(\Ru{(P_{\|\Gamma\|})}\times\Ru{(P_{\|\Gamma\|}^{-})}\bigr)
\leftthreetimes
\bigl((L_{\langle\Gamma\rangle^{\perp}}'T^{|\Gamma|_{\ZZ}}
\times L_{\langle\Gamma\rangle^{\perp}}'T^{|\Gamma|_{\ZZ}})
\cdot\diag L_{|\Gamma|}\bigr)\\
\dim Y&=\dim G-\dim L_{\langle\Gamma\rangle^{\perp}}+\dim\Gamma
\end{align*}
\end{theorem}

\begin{proof}
Take $\gamma\in(\intr\C_Y)\cap\V$. Then
\begin{displaymath}
\gamma(t)=\proj{\sum\gamma(t)v_{\mu}\otimes v_{-\mu}}=\proj{\sum
t^{\langle\gamma,\mu\rangle}v_{\mu}\otimes v_{-\mu}}\to y
\qquad\mbox{as }t\to0
\end{displaymath}
whence $y\in\overline{T}$. Moreover,
$f_{\lambda}(\gamma(t))=
t^{\langle\lambda,\gamma\rangle}c_{\lambda}$,
$c_{\lambda}=f_{\lambda}(e)\ne0$, $\forall\lambda\in\X$.

Choose a closed orbit $Y_0\subseteq\overline{Y}$ corresponding
to a certain vertex $\lambda_0\in C$ of~$\Gamma$. In the
notation of~\S\ref{proj.comp}, $y\in\open{X}$, whence
$y\in Z$. We have to prove that $y\in Y\cap Z$.

It follows from the local structure of $\open{X}$
(\S\ref{loc.struct}) that $Y\cap Z$ is an $(L\times
L)$-orbit in $Z$ and $\C_{Y\cap Z}=\C_Y$ is a face
of~$\Sigma_0^{\vee}$. The ideal of $\overline{Y\cap Z}$ (the closure in~$Z$) as an
$(L\times L)$-submodule in $\CC[Z]$ is given by its highest weight
vectors~$f_{\lambda}$,
$\lambda\in\Sigma_0\setminus\C_Y^{\perp}=
\Sigma_0\setminus|\Gamma|$. Thence $f_{\lambda}(y)=0\iff
f_{\lambda}|_{Y\cap Z}=0$.

Furthermore, for all $(b,b^{-})\in B\times B^{-}$ we have
\begin{displaymath}
(b,b^{-})f_{\lambda}(y)=f_{\lambda}(\proj{b^{-1}e_{\Gamma}b^{-}})=
f_{\lambda}(b_{\Gamma}^{-1}e_{\Gamma}b_{\Gamma}^{-})
\end{displaymath}
(Here $b^{\pm}=u_{\Gamma}^{\pm}b_{\Gamma}^{\pm}$ are the
decompositions in
$B^{\pm}=\Ru{\bigl(P_{|\Gamma|}^{\pm}\bigr)}\leftthreetimes
B_{|\Gamma|}^{\pm}$, $B_{|\Gamma|}^{\pm}=B^{\pm}\cap
L_{|\Gamma|}$.) The latter expression may be regarded as the
value of a function of weight $(-\lambda,\lambda)$ at
$b_{\Gamma}^{-1}b_{\Gamma}^{-}\in L_{|\Gamma|}$, since the
identity maps to $e_{\Gamma}$ under the representation
$L_{|\Gamma|}:V_{\Gamma}$. As the two ``big cells''
$B_{|\Gamma|}B_{|\Gamma|}^{-},B_{|\Gamma|}^{-}B_{|\Gamma|}
\subseteq L_{|\Gamma|}$ intersect in a dense subset, for almost
all $(b,b^{-})$ there exist
$(c_{\Gamma}^{-},c_{\Gamma})\in B_{|\Gamma|}^{-}\times
B_{|\Gamma|}$ such that
$f_{\lambda}(b_{\Gamma}^{-1}e_{\Gamma}b_{\Gamma}^{-})=
f_{\lambda}(c_{\Gamma}^{-}e_{\Gamma}c_{\Gamma})=
d_{\lambda}f_{\lambda}(y)$, $d_{\lambda}\ne0$. (Recall that
$f_{\lambda}$ is an eigenfunction for $B_{|\Gamma|}^{-}\times
B_{|\Gamma|}$.)

It follows that $f_{\lambda}(y)=0$ iff $\CC[Z]_{(\lambda)}$
vanishes at~$y$, whence $y\in\overline{Y\cap Z}$ and is not
contained in smaller orbits. This implies $y\in Y\cap Z$.

Now we pass to stabilizers. For $(g,h)\in G\times G$ we have:
$(g,h)y=y$ iff $ge_{\Gamma}h^{-1}$ is proportional to
$e_{\Gamma}$ iff
\begin{enumerate}
\renewcommand{\theenumi}{\textup{(\arabic{enumi})}}
\renewcommand{\labelenumi}{\theenumi}
\item\label{Im} $gV_{\Gamma}=V_{\Gamma}$,
\item\label{Ker} $hV_{\Gamma}'=V_{\Gamma}'$,
\item\label{scalar} the actions of $g,h$ on $V_{\Gamma}\cong
V/V_{\Gamma}'$ differ by a scalar multiple.
\end{enumerate}
The condition~\ref{Im} means that $g\in P_{\|\Gamma\|}$,
\ref{Ker}~$\iff h\in P_{\|\Gamma\|}^{-}$, and the kernels of the
actions $P_{\|\Gamma\|}^{\pm}:V_{\Gamma}\cong V/V_{\Gamma}'$ are
$\Ru{\bigl(P_{\|\Gamma\|}^{\pm}\bigr)}\leftthreetimes
(L_{\langle\Gamma\rangle^{\perp}}'T^{\langle\Gamma\rangle_{\ZZ}})$,
i.e., we may assume below that $g,h\in L_{|\Gamma|}$. But now
\ref{scalar} just means that, up to multiplying by one and the
same element of~$L_{|\Gamma|}$, $g,h\in T^{|\Gamma|_{\ZZ}}$.

The formula for $\dim Y$ easily follows from the structure of
$(G\times G)_y$ and of~$L_{\|\Gamma\|}'$.
\end{proof}

\begin{Remark}
The theorem generalizes the results of de~Concini--Procesi for
wonderful and regular completions (see \cite{CP1}, \cite{CP2},
\cite{inf(Bruhat)}) and those of Putcha--Renner \cite{irr.mon},
\cite{BN-mon}, \cite{ss.var}, and Vinberg~\cite[Thm.7]{semigr}
for algebraic semigroups. A direct link to algebraic semigroups
is provided by considering the cone over $X$ in $\VV$. It is an
algebraic semigroup, whose group of invertibles is the extension
of $G$ by homotheties (cf.~\cite{ss.var}). The idea of
computing stabilizers is taken from \cite[\S7]{semigr}.

The results close to Proposition~\ref{orb<->faces} and
Theorem~\ref{orb} were obtained for \emph{normal} affine and
projective embeddings of $G$ (and even in a more general context)
by Alexeev and Brion \cite{red.var.aff}, \cite{red.var.pro}. When
the preliminary text of this paper was written, the author knew
about the paper of Kapranov~\cite{hyper} in which all assertions
of Proposition~\ref{orb<->faces} and Theorem~\ref{orb}, except for
stabilizers and dimensions, where proved, see \cite[2.4.2]{hyper}.
However the proof therein seems to be incomplete. The author is
indebted to M.~Brion for this reference.
\end{Remark}

\begin{Corollary}
$\overline{T}$ intersects all $(G\times G)$-orbits in~$X$.
$T$-orbits in the intersection of $\overline{T}$ with a
$(G\times G)$-orbit are permuted by $W$ transitively.
\end{Corollary}
\begin{Remark}
An assertion similar to the first one holds for any spherical variety
\cite[2.4]{var.sph}. A simplest ``transcendental'' proof of the
corollary is obtained by closing in $X$ the Cartan decomposition
$G=KTK$, where $K\subset G$ is a maximal compact subgroup.
One can give an ``algebraic'' proof in the same way by considering
the Iwahori decomposition of $G\bigl(\CC((t))\bigl)$ instead of
the Cartan decomposition, cf.~\cite[2.4, Exemple~2]{var.sph}.

On the other hand, this corollary can in used to obtain a simple
proof of Theorem~\ref{orb} as follows. It is an easy exercise in
toric geometry (cf.~\S\ref{polytope}) that $y=\proj{e_{\Gamma}}$
over \emph{all} faces $\Gamma\subseteq\P$ form a complete set of
$T$-orbit representatives in $\overline{T}$. Thus these $y$
represent all $(G\times G)$-orbits in~$X$. Now it is easy to
deduce from the structure of the respective orbit $Y$ as a
homogeneous space given by $(G\times G)_y$ that $Y^{\diag{T}}$ is
a union of $T$-orbits permuted by $W$ transitively. Therefore
$y=\proj{e_{\Gamma}}$ over those $\Gamma\subseteq\P$ with
$(\intr\Gamma)\cap C\ne\emptyset$ form a complete set of $(G\times
G)$-orbit representatives in~$X$.

However we are also interested in the combinatorial data of the
embedding theory (colored cones) related to these orbits, so our
proof is different.
\end{Remark}

\section{Normality}
\label{norm}

The questions of normality and smoothness are of local nature.
Thus it suffices to examine them at points of some closed orbit
$Y_0\subset X$. A general normality criterion is essentially
contained in Propositions~\ref{BLV-slice},\ref{cone} and the
subsequent remarks.

We say that weights $\mu_1,\dots,\mu_s\in\X^{+}$
\emph{$G$-generate} a semigroup~$S$ if $S$ consists of all
highest weights $k_1\mu_1+\dots+k_s\mu_s-\beta$
($\beta\in\ZZ_{+}\Pi$) of $G$-modules $V(\mu_1)^{\otimes k_1}
\otimes\dots\otimes V(\mu_s)^{\otimes k_s}$,
$k_1,\dots,k_s\in\ZZ_{+}$.

\begin{proposition}\label{crit.norm}
In the notation of~\S\ref{proj.comp}, let
$\mu_1,\dots,\mu_s\in\X^{+}_L$ be all highest weights of $L:V$,
except~$\lambda_0$. Then $X$ is normal at points of~$Y_0$ iff
$\mu_1-\lambda_0,\dots,\mu_s-\lambda_0$ $L$-generate
$\X\cap\Sigma_0$.
\end{proposition}
\begin{proof}
Note that $X$ is normal at points of~$Y_0$ iff $Z$~is normal.
By Proposition~\ref{BLV-slice},
$Z\embeds\VV_0=\bigoplus\LO(V_L(\mu_j-\lambda_0))$, hence
$\CC[Z]\subseteq\CC[L]$ is generated by the components
$\CC[L]_{(\mu_j-\lambda_0)}$. By Proposition~\ref{tails},
$\CC[Z]=\bigoplus_{\lambda\in S}\CC[L]_{(\lambda)}$, where
$S\subseteq\X^{+}_L$ is the semigroup $L$-generated
by~$\mu_j-\lambda_0$, $j=1,\dots,s$. Now the assertion follows
from the corollary of Proposition~\ref{cone}.
\end{proof}

In order to apply this criterion, one requires the information
on decomposing $V$ into simple $L$-modules and on decomposing
tensor products of $L$-modules. The first problem is eliminated
by the following lemma:

\begin{lemma}\label{L-gens}
The weights $\lambda_1-\lambda_0,\dots,\lambda_m-\lambda_0$,
$-\alpha_j\in-(\Pi\setminus\Pi_L)$ are highest weights of
$L:V(-\lambda_0)$ and they $L$-generate all highest weights of
$L:V(-\lambda_0)$.
\end{lemma}
\begin{proof}
Observe that
$v_{\lambda_i},e_{-\alpha_j}v_{\lambda_0}$ are highest weight vectors of
$L:V$, whence the first assertion. Furthermore, for every
$i=0,\dots,m$ there exists a unique $L$-submodule
$V_L(\lambda_i)\embeds V$ (generated by~$v_{\lambda_i}$) and
\begin{displaymath}
V=\sum_{k,i}\underbrace{\g\cdots\g}_k v_{\lambda_i}=
\sum_{k,i}\underbrace{\p^{-}\cdots\p^{-}}_k\cdot V_L(\lambda_i)
\end{displaymath}
The algebra $\p^{-}$ is generated by $\lv$ and~$e_{-\alpha_j}$,
$\alpha_j\notin\Pi_L$, and it contains simple $L$-submodules
$\g_L(-\alpha_j)\cong V_L(-\alpha_j)$ with highest weight vectors
$e_{-\alpha_j}$. Therefore
\begin{displaymath}
V=\sum_{n,i,j_1,\dots,j_n}
\g_L(-\alpha_{j_1})\cdots\g_L(-\alpha_{j_n})\cdot
V_L(\lambda_i)
\end{displaymath}
But $\g_L(-\alpha_{j_1})\cdots\g_L(-\alpha_{j_n})\cdot
V_L(\lambda_i)$ is a quotient module of
$V_L(-\alpha_{j_1})\otimes\dots\otimes V_L(-\alpha_{j_n})\otimes
V_L(\lambda_i)$. This yields the second assertion of the lemma.
\end{proof}

Thus in Proposition~\ref{crit.norm} one may replace
$\mu_1-\lambda_0,\dots,\mu_s-\lambda_0$ by
$\lambda_1-\lambda_0,\dots,\lambda_m-\lambda_0$,
$-\alpha_j\in-(\Pi\setminus\Pi_L)$. (Although sometimes it is
more convenient to operate with all $L$-highest weights if they
are known.) However the problem of decomposing tensor products
of modules with $L$-generating highest weights into simple
$L$-modules remains. (For normality one has to obtain all
generators of $\X\cap\Sigma_0$ among highest weights of all
occurring simple $L$-submodules.) That is why this criterion is
not really effective in the general case. However it implies
simple sufficient conditions.
\begin{Corollary}
$X$ is normal (at points of~$Y_0$) if
$\lambda_1-\lambda_0,\dots,\lambda_m-\lambda_0$ generate
$\X\cap\Sigma_0$.
\end{Corollary}

The necessary condition for normality obtained by Renner for
reductive algebraic semigroups (see~\cite[6.4]{ss.mon}) extends to
projective compactifications.
\begin{proposition}
If $X$ is normal, then $\overline{T}$ is normal.
\end{proposition}
\begin{proof}
Replacing $G$ by~$L$, $X$~by~$Z$, $\overline{T}$~by the
closure of $T$ in~$L$, we reduce the problem to the case of an
affine embedding. We have
$Z\subseteq\VV_0=\bigoplus\LO(V_L(\nu_j))$, where
$\nu_j\in\X^{+}_L$ $L$-generate~$\Sigma_0\cap\X$. We can
increase $\VV_0$ by adding new $L$-highest weights in such a way
that $\nu_j$ will generate $\Sigma_0\cap\X$. Then
$Z=\overline{L}\subseteq\VV_0$ will not change (the highest
weights of $\CC[Z]$ are the same), but now the weights of
$T:V_L(\nu_j)$ will generate the semigroup $W_L\Sigma_0\cap\X$
of all lattice points of~$W_L\Sigma_0$, the cone of
$\P$ at~$\lambda_0$. The respective semigroup algebra is the
coordinate algebra of $\overline{T}\subseteq Z$, whence
$\overline{T}$~is normal.
\end{proof}
\begin{Remark}
This condition can be effectively verified (\S\ref{polytope}),
and in the case, where all vertices of $\P$ are regular weights,
the condition coincides with the general normality criterion.
However this condition is not sufficient in the general case.
\end{Remark}

\begin{Example}
Let $G=\Sp_4$ and the highest weights of $V$ be
$\{\lambda_0,\lambda_1\}=\{3\omega_1,2\omega_2\}$. The weight
polytope $\P$ is given at Fig.~(a), the highest weights are
indicated by bold dots.
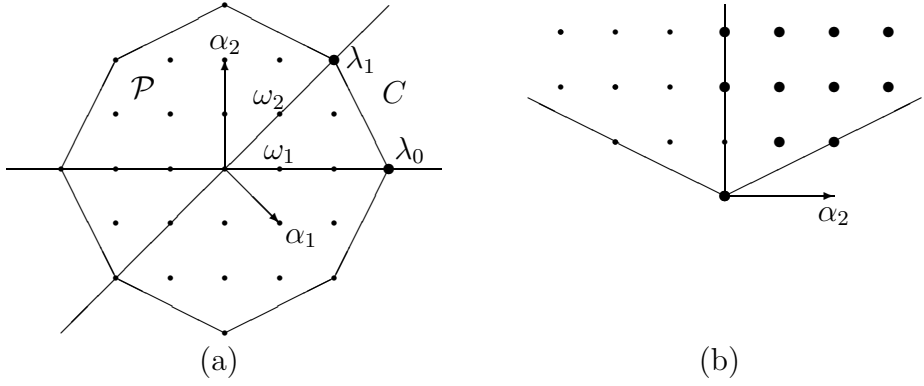
\begin{center}
\begin{tabular}{c}
\unitlength 0.40ex
\linethickness{0.4pt}
\begin{picture}(79.67,60.17)
\put(70.00,60.00){\line(-1,-1){60.00}}
\put(0.00,30.00){\line(1,0){79.67}}
\put(40.00,30.00){\vector(0,1){20.00}}
\put(40.00,30.00){\vector(1,-1){10.00}}
\put(70.00,30.00){\line(-1,2){10.00}}
\put(40.00,60.00){\line(2,-1){20.00}}
\put(10.00,30.00){\line(1,2){10.00}}
\put(40.00,60.00){\line(-2,-1){20.00}}
\put(70.00,30.00){\line(-1,-2){10.00}}
\put(40.00,0.00){\line(2,1){20.00}}
\put(10.00,30.00){\line(1,-2){10.00}}
\put(40.00,0.00){\line(-2,1){20.00}}
\put(40.00,30.00){\circle*{1.00}}
\put(50.00,30.00){\circle*{1.00}}
\put(60.00,30.00){\circle*{1.00}}
\put(30.00,30.00){\circle*{1.00}}
\put(20.00,30.00){\circle*{1.00}}
\put(10.00,30.00){\circle*{1.00}}
\put(40.00,40.00){\circle*{1.00}}
\put(40.00,50.00){\circle*{1.00}}
\put(40.00,60.00){\circle*{1.00}}
\put(50.00,40.00){\circle*{1.00}}
\put(50.00,50.00){\circle*{1.00}}
\put(30.00,40.00){\circle*{1.00}}
\put(30.00,50.00){\circle*{1.00}}
\put(20.00,40.00){\circle*{1.00}}
\put(20.00,50.00){\circle*{1.00}}
\put(40.00,20.00){\circle*{1.00}}
\put(40.00,10.00){\circle*{1.00}}
\put(50.00,20.00){\circle*{1.00}}
\put(50.00,10.00){\circle*{1.00}}
\put(60.00,20.00){\circle*{1.00}}
\put(60.00,10.00){\circle*{1.00}}
\put(30.00,20.00){\circle*{1.00}}
\put(30.00,10.00){\circle*{1.00}}
\put(20.00,20.00){\circle*{1.00}}
\put(20.00,10.00){\circle*{1.00}}
\put(40.00,0.00){\circle*{1.00}}
\put(70.00,30.00){\circle*{2.00}}
\put(60.00,40.00){\circle*{1.00}}
\put(60.00,50.00){\circle*{2.00}}
\put(40.00,51.00){\makebox(0,0)[cb]{$\alpha_2$}}
\put(51.00,19.00){\makebox(0,0)[lt]{$\alpha_1$}}
\put(50.00,31.00){\makebox(0,0)[cb]{$\omega_1$}}
\put(71.00,31.00){\makebox(0,0)[lb]{$\lambda_0$}}
\put(62.00,50.00){\makebox(0,0)[lc]{$\lambda_1$}}
\put(71.00,44.00){\makebox(0,0)[cc]{$C$}}
\put(25.00,45.00){\makebox(0,0)[cc]{$\P$}}
\put(51.00,41.00){\makebox(0,0)[rb]{$\omega_2$}}
\end{picture}
\\ (a)
\end{tabular}
\hfill
\begin{tabular}{c}
\unitlength 0.40ex
\linethickness{0.4pt}
\begin{picture}(72.00,60.17)
\put(6.00,45.00){\circle*{1.00}}
\put(6.00,55.00){\circle*{1.00}}
\put(16.00,35.00){\circle*{1.00}}
\put(16.00,45.00){\circle*{1.00}}
\put(36.00,25.00){\circle*{2.00}}
\put(36.00,35.00){\circle*{1.00}}
\put(26.00,55.00){\circle*{1.00}}
\put(16.00,55.00){\circle*{1.00}}
\put(26.00,45.00){\circle*{1.00}}
\put(16.00,45.00){\circle*{1.00}}
\put(26.00,35.00){\circle*{1.00}}
\put(16.00,35.00){\circle*{1.00}}
\put(6.00,55.00){\circle*{1.00}}
\put(36.00,25.00){\circle*{2.00}}
\put(46.00,35.00){\circle*{2.00}}
\put(56.00,35.00){\circle*{2.00}}
\put(36.00,25.00){\line(0,1){35.00}}
\put(46.00,45.00){\circle*{2.00}}
\put(46.00,55.00){\circle*{2.00}}
\put(56.00,45.00){\circle*{2.00}}
\put(56.00,55.00){\circle*{2.00}}
\put(66.00,45.00){\circle*{2.00}}
\put(66.00,55.00){\circle*{2.00}}
\put(36.00,45.00){\circle*{2.00}}
\put(36.00,55.00){\circle*{2.00}}
\put(36.00,25.00){\vector(1,0){20.00}}
\put(56.00,23.00){\makebox(0,0)[ct]{$\alpha_2$}}
\put(36.00,25.00){\line(2,1){36.00}}
\put(36.00,25.00){\line(-2,1){36.00}}
\end{picture}
\\ (b)
\end{tabular}
\end{center}
Here $L=\SL_2\times\CC^{*}$, $\Delta_L=\{\pm\alpha_2\}$.
The weight semigroup of $\CC[\overline{T}]$ (the closure is
taken in~$Z$) is indicated by dots at Fig.~(b). Bold dots
indicate the subsemigroup of highest weights of~$\CC[Z]$
(which is easy to compute using the Clebsch--Gordan formula).
Now we can see that $\overline{T}$ is normal and $Z$ is not.
\end{Example}

\section{Smoothness}
\label{smooth}

In the theory of spherical varieties, smoothness is usually a
much more subtle property that normality. The general smoothness
criterion \cite[4.2]{geo(sph)} is rather intricate.
Surprisingly, for projective compactifications of reductive
groups, it is easier to verify smoothness, than normality. We
retain the notation of~\S\S\ref{proj.comp},\ref{norm}.

\begin{theorem}\label{reg.pt}
$X$ is smooth at points of~$Y_0$ iff
$L\cong\GL_{n_1}\times\dots\times\GL_{n_p}$ is a \emph{direct}
product, and only polynomial representations of $L$ occur in the
decomposition of the $L$-module $V_0(-\lambda_0)$
(here $V=\CC v_{\lambda_0}\oplus V_0$ as above), and all minimal
representations of factors occur among them.
In the language of roots and weights, this amounts to the
following conditions:
\begin{enumerate}
\renewcommand{\theenumi}{\textup{(\arabic{enumi})}}
\renewcommand{\labelenumi}{\theenumi}
\item\label{Levi-type} All simple components of $L$ are of
type~$\Aa$, and there are no more than $\dim Z(L)$ of
them.
\item\label{simpl} $\Sigma_0$ is simplicial, and moreover, is
generated by a basis of~$\X$.
\item\label{GL} One can enumerate the simple roots in order of
their positions at Dynkin diagrams of connected components 
$\{\alpha^{(k)}_1,\dots,\alpha^{(k)}_{n_k-1}\}$
of~$\Pi_L$, $k=1,\dots,q$, and partition the basis of the free
semigroup $\X\cap\Sigma_0$ into subsets
$\{\pi^{(k)}_1,\dots,\pi^{(k)}_{n_k}\}$, $k=1,\dots,p$, $p\geq
q$, in such a way that
$\langle\pi^{(k)}_j,{\alpha^{(k){\vee}}_j}\rangle=1$,
$\pi^{(k)}_j\perp\Pi_L\setminus\{\alpha^{(k)}_j\}$,
$n_k\pi^{(k)}_j-j\pi^{(k)}_{n_k}\perp
\pi^{(1)}_{n_1},\dots,\pi^{(p)}_{n_p}$, $\forall j,k$.
\item\label{min} Among the weights
$\lambda_1-\lambda_0,\dots,\lambda_m-\lambda_0$,
$-\alpha_j\in-(\Pi\setminus\Pi_L)$ there occur all
$\pi^{(k)}_1$, $1\leq k\leq p$.
\end{enumerate}
\end{theorem}
\begin{proof}
Observe that $X$ is smooth (at points of~$Y_0$) iff $Z$~is
smooth, i.e., the problem is reduced to affine embeddings.

As above, let $\mu_1,\dots,\mu_s$ be the highest weights of
$L:V$. We have an embedding
$Z\embeds\VV_0=\bigoplus\LO(V_L(\mu_j-\lambda_0))$,
$y_0\mapsto0$.

If $Z$ is smooth, then by Luna's fundamental lemma
(\S\ref{fund}), $Z$~projects onto $T_0Z$ isomorphically under an
$(L\times L)$-equivariant projection $\VV_0\to T_0Z$.
Renumeration of~$\mu_j$ yields
\begin{displaymath}
Z\cong\LO(V_L(\mu_1-\lambda_0))\oplus\dots\oplus
\LO(V_L(\mu_p-\lambda_0)),\qquad p\leq s
\end{displaymath}

Let $e\mapsto(e_1,\dots,e_p)$ under this isomorphism. The
projection $Z\to\LO(V_L(\mu_k-\lambda_0))$ maps the dense orbit
$L\subset Z$ onto the dense orbit~$Le_k$, whence
$e_k$ is a nonzero scalar operator. After rescaling
the above isomorphism, we may assume $e_k$ to be the identity
operator. Then the projection maps $L$ homomorphically onto
$\GL(V_L(\mu_k-\lambda_0))$. By a dimension argument,
$L\cong\GL_{n_1}\times\dots\times\GL_{n_p}$, $n_k=\dim
V_L(\mu_k-\lambda_0)$, and $\mu_k-\lambda_0$ are the highest
weights of minimal representations of $\GL_{n_k}$,
$k=1,\dots,p$.

The semigroup $\X\cap\Sigma_0$ is $L$-generated by all
$\mu_k-\lambda_0$ and consists of highest weights of all
polynomial representations of~$L$ (i.e., those extendible
to~$Z$). The description of polynomial representations of
$\GL_n$ implies that $\X\cap\Sigma_0$ is freely generated by the
weights~$\pi^{(k)}_j$, where $\pi^{(k)}_1$ is the highest weight
of a minimal representation of~$\GL_{n_k}$ and
$V_L(\pi^{(k)}_j)= \bigwedge^jV_L(\pi^{(k)}_1)$. Conditions
\ref{Levi-type}--\ref{min} are easily deduced in view of
Lemma~\ref{L-gens}.

Conversely, if conditions \ref{Levi-type}--\ref{min} are
satisfied, then $L\cong\GL_{n_1}\times\dots\times\GL_{n_p}$, and
highest weights in $\CC[Z]$ correspond to all polynomial
representations of~$L$. Therefore
$Z\cong\Mat_{n_1}\times\dots\times\Mat_{n_p}$, i.e., it is
smooth.
\end{proof}
\begin{Remark}
In the case of a regular weight~$\lambda_0$, $L=T$, and the
smoothness criterion is reduced to conditions \ref{simpl}
and~\ref{min}: $X$~is smooth (at points of~$Y_0$) $\iff$
$\Sigma_0$~is generated by the basis
$\pi^{(1)},\dots,\pi^{(p)}$ of~$\X$ and $\lambda_0+\pi^{(k)}$
are in the weight system of $T:V$.
\end{Remark}

\begin{Example}
Let $G=\SO_{2l+1}$ and $V=V(\lambda_0)$ be an irreducible
representation of the fundamental highest weight
$\lambda_0=\omega_i$. For $i<l$,  $X$~will be singular, because
it violates condition~\ref{Levi-type} (or \ref{GL} for $l=2$).

In case of the spinor representation, we have
$\lambda_0=(\eps_1+\dots+\eps_l)/2$, where
$\pm\eps_1,\dots,\pm\eps_l$ are the nonzero
weights of the tautological representation of $\SO_{2l+1}$.
The weights of the spinor representation are
$(\pm\eps_1\pm\dots\pm\eps_l)/2$,
$\Pi=\{\eps_1-\eps_2,\dots
\eps_{l-1}-\eps_l,\eps_l\}$,
$\Pi_L=\Pi\setminus\{\eps_l\}$. The highest weights of
$L:V$ are
$\mu_1=(\eps_1+\dots+\eps_{l-1}-\eps_l)/2$,
$\mu_2=(\eps_1+\dots-\eps_{l-1}-\eps_l)/2$,\dots,
$\mu_l=(-\eps_1-\dots-\eps_{l-1}-\eps_l)/2$.
The vectors $\pi_k=\mu_k-\lambda_0$ generate $\Sigma_0$ and form
a basis of $\X=\langle\eps_1,\dots,\eps_l\rangle$.
Indeed,
$\pi_1=-\eps_1,\pi_2=-\eps_1-\eps_2,\dots,
\pi_l=-\eps_1-\dots-\eps_l$. It is also easy to
see that condition~\ref{GL} is verified. Thus $X$ is smooth.
\end{Example}

\begin{Remark}
Using a more subtle version of the local structure in
the neighborhood of a non-closed orbit, one can obtain
conditions of normality and smoothness at any point of~$X$.
\end{Remark}

\section{Examples}
\label{examples}

Here we illustrate the general theorems proven above by
describing geometric properties of equivariant
compactifications of simple algebraic groups in the spaces of
projective linear operators of fundamental and adjoint
representations. In each case, we describe the structure of the
orbit set and examine normality and smoothness of the
compactification.

Our results are presented in Tables~\ref{Hasse}--\ref{sing}. We
consider a representation of a simple algebraic group $G$,
indicated up to isomorphism in the column ``Group'', in a module
$V=V(\lambda_0)$, where $\lambda_0$ is a fundamental weight or
the highest root. Fundamental representations are denoted in
the column ``Module'' by indicating the highest
weight~$\lambda_0$ and those adjoint representations which are
not fundamental are denoted by the symbol~$\Ad$. The numeration
of simple roots and of fundamental weights, respectively, is
taken from \cite{Lie&alg} (so that $V(\omega_1)$ always has the
minimal dimension).

We consider a $(G\times G)$-variety
$X=\overline{\Ad{G}}\subseteq\PP(\LO(V))$. In the column
``Orbits'' we indicate the dimensions and the Hasse graph of $(G\times
G)$-orbits in~$X$. For classical groups, we also indicate orbit
representatives (given by projectors $V\onto V_{\Gamma}$; if
$V=\g$, then $V_{\Gamma}$ is the center of
$\Ru{(\p_{|\Gamma|})}$ or~$\Ru{(\p_{\|\Gamma\|})}$). These data
can be easily derived from Theorem~\ref{orb} and the following
lemma (cf.~\cite[Thm.\,2]{irr.mon}):
\begin{lemma}
Let $V=V(\lambda_0)$, $\lambda_0\in\X^{+}$, and $L\subseteq
P=P(\lambda_0)$ be the standard Levi subgroup. The faces of
$\P=\P(V)$ intersecting $C$ are of the form
$\Gamma=\P\cap(\lambda_0+\langle\Pi_{\Gamma}\rangle)$, where
$\langle\Pi_{\Gamma}\rangle\subseteq\Pi$ is a subsystem of
simple roots such that no connected component of $\Pi_{\Gamma}$
is contained in~$\Pi_L$, or $\Pi_{\Gamma}=\emptyset$.
Furthermore, $|\Gamma|=\langle\Pi_{\Gamma}\rangle$, the systems
of simple roots of $L_{|\Gamma|}$ and of
$L_{\langle\Gamma\rangle^{\perp}}$ are $\Pi_{\Gamma}$ and
$\Pi_{\Gamma^{\perp}}=\Pi_L\cap(\Pi_{\Gamma})^{\perp}$,
respectively, and $e_{\Gamma}$ is the
$L_{|\Gamma|}$-equivariant projector onto
$V_{L_{|\Gamma|}}(\lambda_0)\subseteq V_G(\lambda_0)$.

In particular, for $\lambda_0=\omega_i$, $\Pi_{\Gamma}$~is a
connected subsystem of~$\Pi$ containing~$\alpha_i$, or
$\Pi_{\Gamma}=\emptyset$, and $\Pi_{\Gamma^{\perp}}=
(\Pi_{\Gamma})^{\perp}\cap\Pi\setminus\{\alpha_i\}$.
\end{lemma}
\begin{proof}
The faces $\Gamma\subseteq\P$, $\intr\Gamma\cap C\ne\emptyset$,
are cut out by supporting functions $\gamma\in-C$. Consider the root
subsystem $\Delta_{\gamma}=\Delta\cap\gamma^{\perp}$ with simple
roots in $\Pi_{\gamma}=\Delta_{\gamma}\cap\Pi$. The space
$|\Gamma|$ is generated by those $\alpha\in\Delta_{\gamma}^{+}$ such
that $e_{-\alpha}v_{\lambda_0}\ne0$, i.e.,
$\alpha\in\Delta_{\gamma}^{+}\setminus\Delta_L^{+}$. In other
words, $\alpha$~has a strictly positive linear expression in a
connected subsystem $\Pi'\subseteq\Pi_{\gamma}$,
$\Pi'\not\subseteq\Pi_L$. All such $\Pi'$ are linearly expressed
in such $\alpha$. Hence $|\Gamma|=\langle\Pi_{\Gamma}\rangle$,
where $\Pi_{\Gamma}$ is the union of connected components of
$\Pi_{\gamma}$ that are not contained in~$\Pi_L$. This yields
the description of~$\Gamma$. The description of $e_{\Gamma}$
stems from $V_{\Gamma}=V_{L_{|\Gamma|}}(\lambda_0)$, and other
assertions of the lemma are evident.
\end{proof}

The Hasse diagrams of orbit sets coincide with those for
irreducible simple algebraic semigroups (which are nothing else,
but the cones over our projective compactifications). The latter
are computed in~\cite{irr.mon}.

We indicate the Levi subgroup of $P=P(\lambda_0)$ in the column
``$L$'' and the nonzero highest weights of $L:V(-\lambda_0)$ in
the column ``$L$-weights''. If there are too many of them, then
we indicate only $-\alpha_j\in-(\Pi\setminus\Pi_L)$ (which
$L$-generate all other highest weights by Lemma~\ref{L-gens}).
We use the following notation: $\eps$~is a fixed basic weight of the
central $1$-torus; $\eps_1,\dots,\eps_n$ are the weights of the
tautological representation in $\CC^n$ of a classical subgroup
of $\GL_n$ ($\eps_{n-i}=-\eps_i$ for the orthogonal and the
symplectic group); $\pi_i=\eps_1+\dots+\eps_i$ is the highest
weight of~$\bigwedge^i\CC^n$; if $L$ is represented as a
quotient of a direct product of several groups, then the weights
of the factors are distinguished by superscripts $'$, $''$, \dots.

In the columns ``Normality'' and ``Smoothness'' we indicate
whether $X$ has the respective property. If normality fails, we
give a reason for it (see
Propositions~\ref{cone},\ref{crit.norm}): if $S$ is the
semigroup of highest weights of $L:\CC[Z]$ and $\mu\notin S$,
but $k\mu\in S$, then we write ``$\notexists\mu$, $\exists
k\mu$''. The following two lemmas are helpful in verifying
normality:

\begin{lemma}
Suppose $G_i:V_i$ ($i=1,\dots,s$) are faithful representations of
connected reductive groups such that $Z(G_i)$ act by
homotheties. Let $Z_i$, $Z$ be the closures of the images
of $G_i$, $G_1\times\dots\times G_s$ in $\LO(V_i)$,
$\LO(V_1\otimes\dots\otimes V_s)$, respectively. Then $Z$ is
normal iff $Z_1,\dots,Z_s$ are normal.
\end{lemma}
\begin{proof}
The algebra $\CC[Z_1\times\dots\times Z_s]$ is multi-graded by
the action of the $s$-dimensional torus
$Z(G_1)\times\dots\times Z(G_s)$, and $\CC[Z]$ is the invariant
algebra of the subtorus $T_0=\{(t_1,\dots,t_s)\mid t_1\dots
t_s=1\}$. Thus normality of $Z$ is implied by normality of
$Z_1,\dots,Z_s$, of direct products and of quotients of normal
varieties \cite[App.\,I, 4.4, II.3.3, Satz~1]{inv}.

On the other hand, the action $Z(G_i):Z_i$ by homotheties lifts
to~$\widetilde{Z_i}$. Let $f_i\in\CC[\widetilde{Z_i}]$ be an
arbitrary homogeneous function. For all $j\ne i$ choose
homogeneous functions $f_j\in\CC[\widetilde{Z_j}]$ of the same
degree and consider $f=f_1\dots
f_s\in\CC[\widetilde{Z_1\times\dots\times Z_s}]^{T_0}$. It is a
rational function on $Z_1\times\dots\times Z_s$ which is
constant on $T_0$-orbits. It is easy to see that all
$T_0$-orbits in $(Z_1\setminus\{0\})\times\dots\times
(Z_s\setminus\{0\})$ are closed in $Z_1\times\dots\times Z_s$,
i.e., generic $T_0$-orbits are closed. Since closed orbits are
separated by invariant polynomials, $f$~is pulled back from a
rational function on~$Z$ which is integral over~$\CC[Z]$
\cite[App.\,I, 3.7, Satz~2, III.3.3, Satz~1]{inv}.
Thus if $Z$ is normal, then $f\in\CC[Z]$, whence
$f_i\in\CC[Z_i]$. Therefore $Z_i$ is normal.
\end{proof}

\begin{lemma}\label{micro}
Let $G$ be one of the groups $\GL_n$,
$\CC^{\times}\times\Sp_{2l}$, $\CC^{\times}\times\SO_{2l}$,
$\CC^{\times}\times\Spin_n$, $\CC^{\times}\times\Ee_l$
($l=6,7$), and $V=V(\mu_0)$, $\mu_0=\pi_k$, $\eps+\omega_1$,
$\eps+\omega_1$, $\eps+\omega_{[n/2]}$, $\eps+\omega_1$,
respectively. Then the closure $Z\subseteq\LO(V)$ of the image
of $G$ is normal.
\end{lemma}
\begin{proof}
Let $S$ be the semigroup of highest weights of $G:\CC[Z]$
(which is $G$-generated by~$\mu_0$). By Proposition~\ref{cone}
and the subsequent remarks, $Z$~is normal iff $S$ is
\emph{saturated}, i.e., $S=\X\cap\Sigma_0$, where $\X=\ZZ{S}$,
and $\Sigma_0=\QQ_{+}S$ is the cone over $C\cap\P$.

In the case $G=\GL_n$, all weights in $S$ are polynomial,
i.e., are of the form $\mu=m_1\eps_1+\dots+m_n\eps_n$,
$m_1\ge\dots\ge m_n\ge0$. It easily follows from the Pieri
formula for the decomposition of $V(\mu)\otimes V(\pi_k)$, where
$\mu$ is any polynomial weight, that $S$ is the saturated
semigroup consisting of all polynomial weights~$\mu$ such that
$m_2+\dots+m_n\ge(k-1)m_1$ and $m_1+\dots+m_n$ is divisible
by~$k$.

In the case $G=\CC^{\times}\times\Sp_{2l}$, the dominant weights
of $V(\mu_0)^{\otimes m}$ are
$\mu=m\eps+m_1\eps_1+\dots+m_l\eps_l$, $m_1\ge\dots\ge m_l\ge0$,
$m-m_1-\dots-m_l\in2\ZZ_{+}$. All of them are in~$S$. Indeed,
$k\eps+\eps_1+\dots+\eps_k$ are highest weights of
$\bigwedge^k V(\mu_0)$, $2\eps$ is a highest weight of
$V(\mu_0)^{\otimes2}$, and all other~$\mu$ are their
$\ZZ_{+}$-linear combinations. Thus $S$ is saturated.

The situation in the case $G=\CC^{\times}\times\SO_{2l}$ is
similar, but dominant weights are of the form
$\mu=m\eps+m_1\eps_1+\dots+m_{l-1}\eps_{l-1}\pm m_l\eps_l$.

In the case $G=\CC^{\times}\times\Spin_{2l+1}$, dominant weights
of $V(\mu_0)^{\otimes m}$ are of the form
$\mu=m\eps+(m_1\eps_1+\dots+m_l\eps_l)/2$, $m\ge m_1\ge\dots\ge
m_l\ge0$, $m\equiv m_1\equiv\dots\equiv m_l\pmod{2}$.
All of them are in~$S$. Indeed, $2\eps+\eps_1+\dots+\eps_k$
($0\le k<l$) are highest weights of $V(\mu_0)^{\otimes2}$,
$\mu_0=\eps+(\eps_1+\dots+\eps_l)/2$, and all other~$\mu$ are
their $\ZZ_{+}$-linear combinations. Thus $S$ is saturated.

In the case $G=\CC^{\times}\times\Spin_{2l}$, dominant weights
of $V(\mu_0)^{\otimes m}$ are of the form
$\mu=m\eps+(m_1\eps_1+\dots+m_l\eps_l)/2$, $m\ge m_1\ge\dots\ge
m_{l-1}\ge|m_l|$, $(l-2)m\ge m_1+\dots+m_{l-1}-m_l$, $m\equiv
m_1\equiv\dots\equiv m_l\pmod{2}$, $lm\equiv
m_1+\dots+m_l\pmod{4}$. These conditions determine a saturated
semigroup, which is generated by
$\mu_0=\eps+\omega_l$, $2\eps+\omega_{l-2k}$ ($1\le k\le l/2$),
$4\eps+\omega_i+\omega_j$ ($0\le i\le j\le l-2$, $i\equiv
j\pmod{2}$), $k\eps+\omega_{l-i}+(k-2)\omega_{l-1}$ ($2\le k\le
i\le l$, $k\equiv i\pmod{2}$),
$(k+2)\eps+\omega_{l-i}+\omega_{l-j}+(k-2)\omega_{l-1}$ ($2\le
k\le i\le j\le l$, $i\equiv j+k\pmod{2}$).
(Here $\omega_0=0$.) It is easy to verify that all of them
are highest weights in the respective
$V(\mu_0)^{\otimes m}$ (where $m$ is the coefficient in~$\eps$).
Thus our semigroup coincides with~$S$.
%

A similar reasoning applies to
$G=\CC^{\times}\times\Ee_l$. One considers a semigroup
$\X\cap\Sigma_0$, where $\Sigma_0=\QQ_{+}(C\cap\P)$. It is
generated by the following weights: $k\eps+\omega_k$,
$(k+9-l)\eps+\omega_k$ ($0\le k\le l-3$, $\omega_0=0$),
$(l-2)\eps+\omega_{l-2}+\omega_l$, $(l-1)\eps+2\omega_{l-2}$,
$4\eps+\omega_{l-2}$, $2\eps+\omega_{l-1}$,
$(l-1)\eps+\omega_{l-1}+2\omega_l$, $l\eps+3\omega_l$,
$3\eps+\omega_l$. One verifies that all of them occur as highest
weights in $V(\mu_0)^{\otimes m}$, whence
$S=\X\cap\Sigma_0$. We omit routine computations.
\end{proof}

\begin{Remark}
Observe that the list of $\mu_0$ restricted to the maximal torus
of $G'$ in Lemma~\ref{micro} is nothing else but the list of
minuscule weights of simple algebraic groups, up to diagram
automorphisms. It is obvious that if (the restriction of)
$\mu_0$ is not minuscule, then the semigroup $S$ is not
saturated, whence $Z$ is not normal. This was first observed by
Faltings \cite{micro}, and he also proved the normality of $Z$
for some minuscule weights. The uniform proof of normality for
any minuscule weight was recently obtained by de Concini
\cite{norm.sgr} using a representation-theoretic lemma from
\cite{proj.norm}. When the preliminary text of this paper was
written, the author learned about the papers \cite{micro},
\cite{norm.sgr}, \cite{proj.norm} from M.~Brion, to whom he is
grateful for these references.
\end{Remark}


{\footnotesize
\begin{longtable}{|c|c|c|}
\caption{Orbits and adherence\label{Hasse}}\\
\hline
Group & Module & Orbits \\
\hline
\endfirsthead
\caption{(continued)}\\
\hline
Group & Module & Orbits \\
\hline
\endhead
\hline
\endfoot

$\Aa_l=\SL_n$ & $\omega_1,\omega_l$ & $Y_r\ni\proj{\phi}$,
$\rk\phi=r$, $\dim Y_r=2nr-r^2-1$ \\*
$n=l+1$ &&
\unitlength 0.70ex
\linethickness{0.4pt}
\begin{picture}(24.00,4.00)
\put(2.00,2.00){\makebox(0,0)[cc]{$Y_n$}}
\put(12.00,2.00){\makebox(0,0)[cc]{${\cdots}$}}
\put(22.00,2.00){\makebox(0,0)[cc]{$Y_1$}}
\put(4.00,2.00){\vector(1,0){5.00}}
\put(15.00,2.00){\vector(1,0){5.00}}
\end{picture} \\
\cline{2-3}

& $\omega_k$ & $Y_{sr}\ni\proj{\bigwedge^k\CC^n\onto
\bigwedge^s\CC^s\otimes\bigwedge^{k-s}\CC^{r-s}}$ \\*
& $1<k<l$ & $\dim Y_{sr}=2nr-r^2-s^2-1$, $\dim
Y_{kk}=2nk-2k^2$ \\*
&&
\unitlength 0.70ex
\linethickness{0.4pt}
\begin{picture}(47.00,24.00)
\put(2.00,2.00){\makebox(0,0)[cc]{$Y_{k-1,n}$}}
\put(15.00,2.00){\makebox(0,0)[cc]{${\cdots}$}}
\put(30.00,2.00){\makebox(0,0)[cc]{$Y_{k-1,k+1}$}}
\put(7.00,2.00){\vector(1,0){5.00}}
\put(18.00,2.00){\vector(1,0){5.00}}
\put(45.00,2.00){\makebox(0,0)[cc]{$Y_{kk}$}}
\put(37.00,2.00){\vector(1,0){5.00}}
\put(2.00,13.00){\makebox(0,0)[cc]{${\vdots}$}}
\put(15.00,12.00){\makebox(0,0)[cc]{${\cdots}$}}
\put(2.00,22.00){\makebox(0,0)[cc]{$Y_{0n}$}}
\put(15.00,22.00){\makebox(0,0)[cc]{${\cdots}$}}
\put(30.00,22.00){\makebox(0,0)[cc]{$Y_{0,k+1}$}}
\put(7.00,22.00){\vector(1,0){5.00}}
\put(18.00,22.00){\vector(1,0){5.00}}
\put(2.00,20.00){\vector(0,-1){6.00}}
\put(2.00,10.00){\vector(0,-1){6.00}}
\put(30.00,13.00){\makebox(0,0)[cc]{${\vdots}$}}
\put(30.00,20.00){\vector(0,-1){6.00}}
\put(30.00,10.00){\vector(0,-1){6.00}}
\put(15.00,20.00){\vector(0,-1){6.00}}
\put(15.00,10.00){\vector(0,-1){6.00}}
\put(7.00,12.00){\vector(1,0){5.00}}
\put(18.00,12.00){\vector(1,0){5.00}}
\end{picture} \\
\cline{2-3}

& $\Ad$ & $\dim Y_0=n^2-1$, $\dim Y_{sr}=n^2-(r-s)^2-2$ \\
\hline
\end{longtable}}

{\footnotesize
\begin{longtable}{|c|c|c|c|c|c|}
\caption{Normality and smoothness\label{sing}}\\
\hline
Group & Module & $L$ & $L$-weights & Normality &
Smoothness \\
\hline
\endfirsthead
\caption{(continued)}\\
\hline
Group & Module & $L$ & $L$-weights & Normality &
Smoothness \\
\hline
\endhead
\hline
\endfoot

$\Aa_l=\SL_n$ & $\omega_1,\omega_l$ & $\GL_l$ & $\pi'_1$ & Yes &
Yes \\
\cline{2-6}

$n=l+1$ & $\omega_k$ & $\PG(\LO_k\times\LO_{n-k})$ &
$\pi'_i+\pi''_i$ & Yes & No \\*
& $1<k<l$ && $i\le k,\ n-k$ && \\
\cline{2-6}

& $\Ad$ & $\CC^{\times}\times\GL_{n-2}$ &
$\eps',2\eps',\pi''_1,\eps'+\pi''_1,$ & Yes & $n\le3$ \\*
&&& $\eps'+\pi''_{n-3}-\pi''_{n-2},$ && \\*
&&& $2\eps'+\pi''_{n-3}-\pi''_{n-2},$ && \\*
&&& $\eps'+\pi''_1+\pi''_{n-3}-\pi''_{n-2}$ && \\
\hline

$\Bb_l=\SO_n$ & $\omega_k$ & $\GL_k\times\SO_{n-2k}$ &
$\pi'_i+\pi'_j+\pi''_m$ & No & No \\
\hline
\end{longtable}}


\end{document}